\documentclass[titlepage,draft,12pt]{article} 
\usepackage{amssymb,amsthm} 
\usepackage[a4paper]{geometry}
\usepackage[usenames,dvipsnames]{xcolor}

\date{}

\newcommand{\ep}{\varepsilon}
\renewcommand{\qed}{{\penalty 10000\mbox{$\quad\Box$}}}
\newcommand{\re}{\mathbb{R}}
\newcommand{\C}{\mathbb{C}}
\newcommand{\n}{\mathbb{N}}

\newcommand{\Ghat}{\widehat{G}}
\newcommand{\Ehat}{\widehat{E}}

\newcommand{\calh}{\mathcal{H}}
\newcommand{\calhr}{\mathcal{H}_{r}}

\newtheorem{thm}{Theorem}[section]
\newtheorem{thmbibl}{Theorem}

\newtheorem{rmk}[thm]{Remark}
\newtheorem{prop}[thm]{Proposition}
\newtheorem{defn}[thm]{Definition}

\newtheorem{lemma}[thm]{Lemma}

\title{Finding the exact decay rate of all solutions to some second 
order evolution equations with dissipation}

\author{Marina Ghisi\vspace{1ex}\\ 
{\normalsize Universit\`a degli Studi di Pisa} \\
{\normalsize Dipartimento di Matematica}\\ 
{\normalsize PISA (Italy)}\\
{\normalsize e-mail: \texttt{ghisi@dm.unipi.it}}
\and
Massimo Gobbino\vspace{1ex}\\ 
{\normalsize Universit\`a degli Studi di Pisa} \\
{\normalsize Dipartimento di Matematica}\\ 
{\normalsize PISA (Italy)}\\  
{\normalsize e-mail: \texttt{m.gobbino@dma.unipi.it}}
\and
Alain Haraux\vspace{1ex}\\ 
{\normalsize Universit\'{e} Pierre et Marie Curie} \\
{\normalsize Laboratoire Jacques-Louis Lions}\\ 
{\normalsize PARIS (France)}\\  
{\normalsize e-mail: \texttt{haraux@ann.jussieu.fr}}}

\begin{document}
\maketitle
\begin{abstract}

	We consider an abstract second order evolution equation with
	damping.  The ``elastic'' term is represented by a self-adjoint
	nonnegative operator $A$ with discrete spectrum, and the nonlinear
	term has order greater than one at the origin.  We investigate the
	asymptotic behavior of solutions.
	
	We prove the coexistence of slow solutions and fast solutions.
	Slow solutions live close to the kernel of $A$, and decay as
	negative powers of $t$ as solutions of the first order equation
	obtained by neglecting the operator $A$ and the second order
	time-derivatives in the original equation.  Fast solutions live
	close to the range of $A$ and decay exponentially as solutions of
	the linear homogeneous equation obtained by neglecting the
	nonlinear terms in the original equation.
		
	The abstract results apply to semilinear dissipative hyperbolic
	equations.
	
\vspace{6ex}

\noindent{\bf Mathematics Subject Classification 2010 (MSC2010):}
35L71, 35L90, 35B40.


\vspace{6ex}

\end{abstract}

 
\section{Introduction}

In this paper we study the precise asymptotic behavior of decaying solutions
to the second order evolution equation
\begin{equation}
	u''(t)+2\delta u'(t)+Au(t)=f(u(t))
	\quad\quad
	\forall t\geq 0,
	\label{pbm:hyp}
\end{equation}
where $\delta>0$ is a real parameter, $A$ is a self-adjoint linear
operator on a Hilbert space $H$, and $f$ is a nonlinear term.

We assume that $A$ is non-negative, but not necessarily strictly
positive, and that its spectrum is a finite or countable set of
eigenvalues without finite accumulation points.  We also assume that
the nonlinear term has order greater than one in the origin, in the
sense that it satisfies inequalities such as
$$|f(u)|\leq K_{0}\left(|u|^{1+p}+|A^{1/2}u|^{1+q}\right)$$
for some positive exponents $p$ and $q$.

As model examples, have in mind such semilinear hyperbolic
equations  as
\begin{equation}
	u_{tt}+u_{t}-\Delta u+|u|^{p}u=0
	\label{eqn:n}
\end{equation}
with Neumann boundary conditions in a bounded domain 
$\Omega\subseteq\re^{n}$, or
\begin{equation}
	u_{tt}+u_{t}-\Delta u-\lambda_{1}(\Omega)u+|u|^{p}u=0
	\label{eqn:d}
\end{equation}
with Dirichlet boundary conditions in a bounded domain
$\Omega\subseteq\re^{n}$, where $\lambda_{1}(\Omega)$ denotes the
first eigenvalue of the Dirichlet Laplacian.  We point out that in
both cases the operator associated to the linear part has a nontrivial
kernel.

This paper is the final step of a project started
with~\cite{ggh:sol-lentes} and~\cite{ggh:casc-par}.
In~\cite{ggh:casc-par} we investigated the corresponding first order
equation
\begin{equation}
	u'(t)+Au(t)=f(u(t))
	\label{pbm:par}
\end{equation}
under analogous assumptions.  The main results obtained in the first
order case are the following.
\begin{enumerate}
	\renewcommand{\labelenumi}{(\arabic{enumi})}
	
	\item \emph{Slow-fast alternative}.  All non-zero solutions to
	(\ref{pbm:par}) which decay to 0 are either slow solutions, in the
	sense that they decay at most as $t^{-1/p}$, or fast solutions
	decaying to 0 exponentially.

	\item \emph{Asymptotic profile of fast solutions}. Every fast 
	solution $u(t)$ is asymptotic to a solution $v(t)$ to the corresponding 
	homogeneous equation
	\begin{equation}
		v'(t)+Av(t)=0
		\label{eqn:par-h}
	\end{equation}
	in the sense that the difference $u(t)-v(t)$ decays faster than both $u$ and $v$. More 
	precisely, we can always take $v(t)$ to be a ``pure'' solution to 
	(\ref{eqn:par-h}) of the form $v(t):=v_{0}e^{-\lambda t}$, where 
	$\lambda$ is an eigenvalue of $A$ and $v_{0}$ is a corresponding 
	eigenvector.
	
	\item \emph{Existence of an open set of slow solutions}.  If
	$\ker(A)$ is nontrivial, and $f$ satisfies a natural sign
	condition, then there exists a nonempty open set of initial data
	giving rise to slow solutions. These solutions live close to the 
	kernel of $A$, and decay as solutions to the ordinary 
	differential equation 
	\begin{equation}
		u'(t)=-|u(t)|^{p}u(t).
		\label{eqn:ode-par}
	\end{equation}

	\item \emph{Existence of families of fast solutions}.  For every
	small enough ``pure'' solution $v(t)=v_{0}e^{-\lambda t}$ to the
	homogeneous equation (\ref{eqn:par-h}), there exists a family of
	fast solutions to (\ref{pbm:par}) asymptotic to $v(t)$.  This
	family has the same structure as the family of solutions to
	(\ref{eqn:par-h}) which are asymptotic to the given pure solution
	$v(t)$.

\end{enumerate}

In the recent paper~\cite{ggh:sol-lentes} we proved the existence of a
nonempty open set of slow solutions for the second order equation
(\ref{pbm:hyp}) under the additional assumption that $f(u)=-\nabla
F(u)$ for a suitable nonnegative functional $F(u)$.  Once again, these
solutions live close to the kernel of $A$ and decay as the solutions to
the first order ordinary differential equation (\ref{eqn:ode-par}).
Roughly speaking, this means that in the slow regime both  operator
$A$ and second order time-derivative can be neglected in
(\ref{pbm:hyp}).  This result extends point~(3) above from the first
order equation (\ref{pbm:par}) to a large class of second order
equations (\ref{pbm:hyp}).

In this paper we extend points (1), (2) and (4).  In
Theorem~\ref{thm:main-alternative} we prove the slow-fast alternative
and we describe the asymptotic profile of fast solutions, which now
behave as solutions to the linear homogeneous equation obtained from
(\ref{pbm:hyp}) by neglecting the nonlinear term.  In
Theorem~\ref{thm:main-exponential} we construct families of fast
solutions with a given asymptotic profile.  The main difference
from~\cite{ggh:casc-par} is that \emph{in the first order case fast
solutions can have infinitely many exponential decay rates},
corresponding to eigenvalues of $A$, while here \emph{in the second
order case only finitely many exponential decay rates can occur}, even
if $A$ has infinitely many distinct eigenvalues.

In the main results of this paper we never require a gradient
structure on the nonlinear term.  On the contrary, our assumptions do
not even guarantee the existence of global solutions for all initial
data in a neighborhood of the origin, and hence our abstract results
apply also to those equations which exhibit coexistence of decaying
solutions and solutions that do not globally exist.  A typical example
is provided by partial differential equations such as (\ref{eqn:n}) or
(\ref{eqn:d}), but with the minus sign in front of the nonlinear term.

Concerning the technique, the literature seems to reveal a shortage of
appropriate tools to tackle questions of this type, even in finite
dimensions.  For example, the classical linearization results \`{a} la
Hartmann (see \cite{hartmann-1,hartmann-2,hartmann-3}) provide a good
description of the dynamic in a neighborhood of a stationary point
(the origin in this case).  On the other hand, the linearization is
realized through homeomorphisms which are just H\"{o}lder continuous,
and therefore they do not preserve decay rates.  More important,
almost all these results seems to deal with the case where eigenvalues
of the linearized equation have real part different from zero, while
we know that a nontrivial kernel is exactly what produces slow
solutions.  Finally, these tools seem to require assumptions on the
dynamics in a whole neighborhood of the origin.

For the same reasons, also geometric tools such as stable or center
manifolds are unlikely to be helpful in answering questions~(1)
through~(4) above.  On the contrary, it is our classification of decay
rates which seems to lead to a better description of the stable
manifold $\mathcal{S}$, when it exists, since what we actually provide
is a ``stratification'' of $\mathcal{S}$ into submanifolds
$\mathcal{S}_{1}\supset
\mathcal{S}_{2}\supset\mathcal{S}_{3}\supset\ldots$ corresponding to
different decay rates.

A first step in the classification of decay rates was done by the 
third author in~\cite{h:ode}, where the full answer is given in the case 
of the scalar ordinary differential equation
\begin{equation}
	u''(t)+u'(t)+|u(t)|^{p}u(t)=0.
	\label{eqn:ode-h}
\end{equation}

The technique exploited in~\cite{h:ode} seems to be
specific to the scalar case .  Nevertheless, this was enough to show
the existence of both slow and fast solutions to the Neumann problem
for (\ref{eqn:n}), since it is enough to consider spatially
homogeneous solutions depending only on $t$, for which (\ref{eqn:n})
reduces to (\ref{eqn:ode-h}).  It was later shown in~\cite{h:hyp} that
the slow-fast alternative holds true for (\ref{eqn:n}), and that the
set of initial data producing fast solutions is closed with empty
interior.  Unfortunately, the proofs of these results seem to exploit
in an essential way the fact that the kernel of the linear part (in
this case the set of constant functions) is an invariant space also
for the nonlinear equation.  Without this assumption, which fails for
example in the Dirichlet case, even the existence of a single slow
solution was open until~\cite{ggh:sol-lentes}.

For all these reasons, in this paper we follow a different path, close
to what we did in the first order case.  We exploit two main tools,
one linear and one nonlinear.  The first one (see
section~\ref{sec:lin}) is a sharp analysis of the linear homogeneous
equation obtained by replacing the right-hand side of (\ref{pbm:hyp})
by a forcing term $g(t)$.  The second one are what we call modified
Dirichlet quotients (see section~\ref{sec:energies}).  They have been
developed in completely different contexts, for example backward
uniqueness results for parabolic differential equations, but they
proved to be fundamental in~\cite{ggh:casc-par} in order to establish
the slow-fast alternative.  Here we need a ``hyperbolic version'' of
Dirichlet quotients, analogous to the quotients introduced
in~\cite{ggh:sol-lentes} and previously
in~\cite{ghisi:JDE2006,gg:JDE2008}.

When we apply the results of~\cite{ggh:sol-lentes} and of the present
paper to the model examples we started with, we end up with a
satisfactory description of the asymptotic behavior of solutions.
More important, this description descends from an abstract framework
which applies in the same way to both ordinary and partial
differential equations, both Neumann and Dirichlet boundary
conditions, both equations with the right and with the wrong sign.

Just as a further example of this flexibility, we present in the last
part of this paper a simple application to a system of ordinary
differential equations describing the so-called finite modes of a
degenerate hyperbolic equation of Kirchhoff type, actually a 
quasi-linear equation.

This paper is organized as follows.  In
Section~\ref{sec:preliminaries} we fix the notation and we introduce
the terminology needed when dealing with decay rates of solutions to
the linear homogeneous equation of order two.  In
section~\ref{sec:statements} we state our main abstract results and we
comment on them.  In Section~\ref{sec:proofs} we develop our tools and
we prove the main results.  In Section~\ref{sec:applications} we
present some simple applications of the abstract theory.

\setcounter{equation}{0}
\section{Notation and preliminaries}\label{sec:preliminaries}

Throughout this paper $H$ denotes a separable Hilbert space, $|x|$ denotes the
norm of an element $x\in H$, and $\langle x,y\rangle$ denotes the
scalar product of two elements $x$ and $y$ in $H$.  We consider a
self-adjoint linear operator $A$ on $H$ with dense domain $D(A)$.  We
assume that $A$ is nonnegative, namely $\langle Au, u\rangle\geq 0$
for every $u \in D(A)$, so that for every $\alpha\geq 0$ the power
$A^{\alpha}u$ is defined provided that $u$ lies in a suitable domain
$D(A^{\alpha})$, which is itself a separable Hilbert space with norm
$$|u|_{D(A^{\alpha})} := \left(|u|^{2} +
|A^{\alpha}u|^{2}\right)^{1/2}.$$

\subsection{The notion of solution}

Let us spend a few words on the notion of solutions to semilinear
equations and inequalities.  Let us start by considering the
linear equation
\begin{equation}
	u''(t)+2\delta u'(t)+Au(t)=g(t),
	\label{eqn:lin-nh}
\end{equation}
with initial data
\begin{equation}
	u(0)=u_{0},
	\hspace{3em}
	u'(0)=u_{1}.
	\label{eqn:u-data}
\end{equation}

There are several ways to introduce a notion of weak solution to 
evolution problems, for example as uniform limits of strong 
solutions, or through integral formulations, but fortunately all of 
them coincide in the case of linear equations such as 
(\ref{eqn:lin-nh}) with initial data in the so-called energy space 
$D(A^{1/2})\times H$.

For our purposes we limit ourselves to forcing terms $g(t)$ which are
defined for every $t\geq 0$ and continuous with values in $H$.  We
recall the following classical result just to state precisely the
regularity of solutions and energy functions.

\begin{thmbibl}[Linear equation -- Existence]\label{thmbibl:linear}
	Let $H$ be a separable Hilbert space, and let $A$ be a self-adjoint
	nonnegative operator on $H$ with dense domain $D(A)$.  Let us
	assume that $g\in C^{0}([0,+\infty),H)$ and $(u_{0},u_{1})\in
	D(A^{1/2})\times H$.
	
	Then problem (\ref{eqn:lin-nh})--(\ref{eqn:u-data}) has a unique
	(weak) solution
	\begin{equation}
		u\in C^{0}\left([0,+\infty),D(A^{1/2})\right)\cap
		C^{1}\left([0,+\infty),H\right).
		\label{th:u-reg}
	\end{equation}
	
	Moreover, the function 
	\begin{equation}
		E(t):=|u'(t)|^{2}+|A^{1/2}u(t)|^{2}
		\label{defn:E}
	\end{equation}
	is of class $C^{1}$ in $[0,+\infty)$, and 
	\begin{equation}
		E'(t)=-4\delta|u'(t)|^{2}+2\langle u'(t),g(t)\rangle
		\quad\quad 
		\forall t\geq 0.
		\label{th:E'}
	\end{equation}
\end{thmbibl}

Next step is considering semilinear equations of the form
\begin{equation}
	u''(t)+2\delta u'(t)+Au(t)=f(u(t)).
	\label{eqn:semilin}
\end{equation}

We assume that $f:B_{R_{0}}\to H$, where $R_{0}>0$ and
\begin{equation}
	B_{R_{0}}:=\left\{u\in D(A^{1/2}):
	|u|_{D(A^{1/2})}<R_{0}\right\}.
	\label{defn:ball}
\end{equation}

We also assume that $f$ is continuous in $B_{R_{0}}$ with respect to
the norm of
$D(A^{1/2})$. We say that $u(t)$ is a solution to 
(\ref{eqn:semilin}) for $t\geq 0$ if there exists 
$g\in C^{0}([0,+\infty),H)$ such that $u(t)$ is a solution to 
(\ref{eqn:lin-nh}) in $[0,+\infty)$, and in addition
$$g(t)=f(u(t))
\quad\quad
\forall t\geq 0.$$

Finally, we consider differential inequalities such as
\begin{equation}
	\left|u''(t)+2\delta u'(t)+Au(t)\strut\right|\leq 
	K_{0}\left(|u(t)|^{1+p}+|A^{1/2}u(t)|^{1+q}\right).
	\label{diff-ineq}
\end{equation}

The notion of solution to these inequalities can be introduced in 
analogy with the case of equation (\ref{eqn:semilin}), as follows. 

\begin{defn}\label{defn:sol}
	\begin{em}
		We say that the function $u(t)$ is a global solution to the
		differential inequality (\ref{diff-ineq}) if
		\begin{itemize}
			\item  $u(t)$ has the regularity stated in 
			(\ref{th:u-reg}),
		
			\item there exists $g\in C^{0}([0,+\infty),H)$ such that
			$u(t)$ is a solution to (\ref{eqn:lin-nh}) in
			$[0,+\infty)$,
		
			\item  the forcing term $g(t)$ satisfies
			\begin{equation}
				|g(t)|\leq
				K_{0}\left(|u(t)|^{1+p}+|A^{1/2}u(t)|^{1+q}\right)
				\quad\quad \forall t\geq 0.
				\label{hp:main-g}
			\end{equation}
		\end{itemize}
	\end{em}
\end{defn}

In particular, the energy (\ref{defn:E}) of all solutions to
(\ref{diff-ineq}) is of class $C^{1}$, and its time-derivative is
given by (\ref{th:E'}).

\subsection{Decay rates for the linear homogeneous equation}\label{sec:decay}

In the sequel we assume that the spectrum $\sigma(A)$ of $A$ is a
finite or countable set of eigenvalues without finite
accumulation points.  Under this assumption, the space $H$ admits a
finite or countable orthonormal system made by eigenvectors of $A$.
We denote this system by $\{e_{k}\}$, where $k$ ranges over some
finite or countable set of indices $\mathcal{K}$.  The corresponding
eigenvalues of $A$ are denoted by $\{\lambda_{k}\}$, so that
$$Ae_{k}=\lambda_{k}e_{k} \quad\quad
\forall k\in\mathcal{K}.$$

We never assume that eigenvalues of $A$ are simple or with finite
multiplicity, so that $\lambda_{k}$'s are not necessarily distinct,
and it could even happen that $\lambda_{k}$ is the same element of
$\sigma(A)$ for infinitely many indices $k$.

For every $\lambda\in\sigma(A)$ we consider the polynomial
\begin{equation}
	z^{2}-2\delta z+\lambda.
	\label{char-pol}
\end{equation}

Then we define the set
$$\mathcal{D}:=\left\{\Re(z):\mbox{$z\in\C$ is a root of 
(\ref{char-pol}) for some $\lambda\in\sigma(A)$}\right\},$$
where $\Re(z)$ denotes the real part of $z$. With some standard 
algebra it is possible to list the elements of $\mathcal{D}$. Indeed, 
for all eigenvalues $\lambda<\delta^{2}$ the  
polynomial (\ref{char-pol}) has two distinct real roots,
which thus provide two elements of $\mathcal{D}$, one less than 
$\delta$ and one greater than $\delta$. All eigenvalues 
$\lambda\geq\delta^{2}$ produce the same element $\delta$ of 
$\mathcal{D}$, as a real root of multiplicity 2 if  
$\lambda=\delta^{2}$, and as the real part of the two distinct 
complex conjugate roots if $\lambda>\delta^{2}$. Since we assumed 
that $\sigma(A)$ has no finite accumulation point, from this list it 
follows that $\mathcal{D}$ is a \emph{finite set}.

The set $\mathcal{D}$ is strongly related to decay rates of solutions 
to the homogeneous linear equation
\begin{equation}
	u''(t)+2\delta u'(t)+Au(t)=0.
	\label{eqn:lin-h}
\end{equation}

In order to make the relation more explicit, we start by considering
\emph{simple modes}, namely solutions to (\ref{eqn:lin-h}) of the form
$u_{k}(t)e_{k}$, where $e_{k}$ is one of the elements of the
orthogonal system, and $u_{k}(t)$ is a solution to the ordinary
differential equation
\begin{equation}
	u_{k}''(t)+2\delta u_{k}(t)+\lambda_{k}u_{k}(t)=0.
	\label{eqn:ode-lin-h}
\end{equation}

The form of the solutions to (\ref{eqn:ode-lin-h}) depends on the
relative order of $\lambda_{k}$ and $\delta^{2}$, as follows.  Let
$(u_{0k},u_{1k})$ be the initial data.
\begin{itemize}
	\item  If $\lambda_{k}>\delta^{2}$, the solution to 
	(\ref{eqn:ode-lin-h}) is
	\begin{equation}
		u_{k}(t)=e^{-\delta t}\left( u_{0k}\cos(\phi_{k}t)+
		\frac{u_{1k}+\delta u_{0k}}{\phi_{k}}
		\sin\left(\phi_{k}t\right)\right),
		\label{sol:uk-sup}
	\end{equation}
	where $\phi_{k}$ denotes the imaginary part of the roots 
	of the polynomial (\ref{char-pol}) with 
	$\lambda=\lambda_{k}$. It always decays as $e^{-\delta 
	t}$, and indeed in this case $\delta$ is the element of 
	$\mathcal{D}$ corresponding to the eigenvalue $\lambda_{k}$ of $A$.

	\item  If $\lambda_{k}=\delta^{2}$, the solution to 
	(\ref{eqn:ode-lin-h}) is
	\begin{equation}
		u_{k}(t)=e^{-\delta t}\left( u_{0k}+(u_{1k}+\delta
		u_{0k})t\right).
		\label{sol:uk-=}
	\end{equation}
	
	Also in this case the exponential term in the decay rate is
	$e^{-\delta t}$, and $\delta$ is the element of $\mathcal{D}$
	corresponding to $\lambda_{k}$.

	\item  If $\lambda_{k}<\delta^{2}$, the solution to 
	(\ref{eqn:ode-lin-h}) is
	\begin{equation}
		u_{k}(t)=
		\frac{u_{1k}+r_{2,k}u_{0k}}{r_{2,k}-r_{1,k}}\cdot 
		e^{-r_{1,k}t}-
		\frac{u_{1k}+r_{1,k}u_{0k}}{r_{2,k}-r_{1,k}}\cdot e^{-r_{2,k}t},
		\label{sol:uk-sub}
	\end{equation}
	where
	\begin{equation}
		r_{1,k}:=\delta-\sqrt{\delta^{2}-\lambda_{k}},
		\hspace{4em}
		r_{2,k}:=\delta+\sqrt{\delta^{2}-\lambda_{k}}.
		\label{defn:r12k}
	\end{equation}
	
	In this case the simple mode is actually the sum of two simple
	modes with different decay rates, described by $r_{1,k}$ and
	$r_{2,k}$, namely by the elements of $\mathcal{D}$ corresponding
	to $\lambda_{k}$ in this range.

\end{itemize}

Therefore, if we limit ourselves to simple modes, all possible decay 
rates involve an exponential term of the form $e^{-rt}$ for some 
$r\in\mathcal{D}$.

The general case is not so different.  Indeed, any solution to
(\ref{eqn:lin-h}) is the sum (or the series) of simple modes.  If we
group together all terms with the same exponential term, we obtain a
decomposition of the form
\begin{equation}
	u(t)=\sum_{r\in D}u_{r}(t).
	\label{defn:ur}
\end{equation}

We point out that this decomposition is unique and involves only a
\emph{finite number} of terms (some of which might be zero).  Each
term is itself a solution to (\ref{eqn:lin-h}), and it can be the
series of countably many simple modes, all with the same exponential
factor.  We specify that terms of type $te^{-\delta t}$, which might
come from (\ref{sol:uk-=}), are grouped together with all other terms
involving $e^{-\delta t}$.

The decomposition (\ref{defn:ur}) motivates the following definitions.

\begin{defn}
	\begin{em}
		Let $r_{0}\in\mathcal{D}$, and let $u(t)$ be a solution to the
		homogeneous equation (\ref{eqn:lin-h}).  Let $u_{r}(t)$ be the
		components of $u(t)$ in the decomposition (\ref{defn:ur}).  We
		say that
		\begin{itemize}
			\item $u(t)$ is \emph{$r_{0}$-fast} if $u_{r}(t)\equiv 0$
			for every $r\leq r_{0}$,
		
			\item $u(t)$ is \emph{$r_{0}$-slow} if $u_{r}(t)\equiv 0$
			for every $r>r_{0}$,
		
			\item $u(t)$ is \emph{$r_{0}$-pure} if $u_{r}(t)\equiv 0$
			for every $r\neq r_{0}$.
		\end{itemize} 
	\end{em}
\end{defn}

\begin{defn}
	\begin{em}
		A pair $(u_{0}.u_{1})\in D(A^{1/2})\times H$ is called
		$r_{0}$-fast (respectively, $r_{0}$-slow or $r_{0}$-pure) if
		the solution to the homogeneous equation (\ref{eqn:lin-h})
		with initial data $(u_{0},u_{1})$ is $r_{0}$-fast
		(respectively, $r_{0}$-slow or $r_{0}$-pure).
	\end{em}
\end{defn}

\begin{rmk}
	\begin{em}
		One could equivalently say that a solution is $r_{0}$-fast if
		it is the sum (or the series) of simple modes with exponential
		factors $e^{rt}$ with $r>r_{0}$.  Analogously, a solution is
		$r_{0}$-slow if it is the sum (or the series) of simple modes
		with exponential factors $e^{rt}$ with $r\leq r_{0}$. We 
		point out that non-zero simple modes with $r=r_{0}$ are 
		allowed in $r_{0}$-slow solutions, but not in $r_{0}$-fast 
		solutions.
	\end{em}
\end{rmk}

\begin{rmk}\label{rmk:lin-h-est}
	\begin{em}
		Let $u(t)$ be a solution to the homogeneous equation 
		(\ref{eqn:lin-h}). Let us assume that $u(t)$ is $r_{0}$-pure 
		or $r_{0}$-fast for some $r_{0}\in\mathcal{D}$. Then it turns 
		out that
		$$\lim_{t\to+\infty}
		\left(|u'(t)|+|u(t)|_{D(A^{1/2})}\right)e^{\gamma t}=0
		\quad\quad
		\forall\gamma<r_{0}.$$
		
		Conversely, if $u(t)$ is $r_{0}$-slow it turns out that
		$$\liminf_{t\to+\infty}
		\left(|u'(t)|+|u(t)|_{D(A^{1/2})}\right)e^{r_{0} t}>0,$$
		unless $u(t)$ is identically 0.
		These results follow in a standard way from the explicit
		expressions (\ref{sol:uk-sup}), (\ref{sol:uk-=}), and
		(\ref{sol:uk-sub}) for the components of $u(t)$.
	\end{em}
\end{rmk}

\subsection{Second order equations as first order systems}\label{sec:defn-system}

In this section we present an alternative description of the
set $\mathcal{D}$ and of the decomposition (\ref{defn:ur}).
Setting $U(t):=(u(t),u'(t))$, the second order equation
(\ref{eqn:lin-h}) can be written as a first order system
\begin{equation}
	U'(t)+\mathcal{A}U(t)=0
	\label{system}
\end{equation}
in the product space $\mathcal{H}:=D(A^{1/2})\times H$, where 
the operator $ \mathcal{A}$ is defined by
\begin{equation}
	\mathcal{A}:=\left(
	\begin{array}{cc}
		0 & -I  \\
		A & 2\delta I
	\end{array}
	\right)			
	\label{defn:cal-A}
\end{equation}
(here $I$ denotes the identity on $H$). It can be proved that 
the spectrum of $\mathcal{A}$ is the union of the roots of 
the polynomials (\ref{char-pol}) when $\lambda$ ranges over
$\sigma(A)$. Therefore, the set 
$\mathcal{D}$ is just the set of real parts of eigenvalues of 
$\mathcal{A}$.

Every element $e_{k}$ of the orthonormal system in $H$ gives
rise to either an $\mathcal{A}$-invariant subspace
$\mathcal{H}_{k}$ of $\mathcal{H}$ of dimension two, or two
$\mathcal{A}$-invariant subspaces $\mathcal{H}_{1,k}$ and
$\mathcal{H}_{2,k}$ of $\mathcal{H}$ of dimension one,
depending on the corresponding eigenvalue $\lambda_{k}$.
\begin{itemize}
	\item If $\lambda_{k}>\delta^{2}$, then $\mathcal{H}_{k}$
	is the two-dimensional subspace generated by $(e_{k},0)$
	and $(0,e_{k})$, where the action of $\mathcal{A}$ has 
	canonical form represented by the matrix
	$$\left(
	\begin{array}{cc}
		\delta & \phi_{k}  \\
		-\phi_{k} & \delta
	\end{array}
	\right),$$
	where $\phi_{k}$ is the same as in (\ref{sol:uk-sup}).

	\item If $\lambda_{k}=\delta^{2}$, then $\mathcal{H}_{k}$
	is the two-dimensional subspace generated by $(e_{k},0)$
	and $(0,e_{k})$, where the action of $\mathcal{A}$ has 
	canonical form represented by the matrix
	$$\left(
	\begin{array}{cc}
		\delta & 1  \\
		0 & \delta
	\end{array}
	\right).$$

	\item If $\lambda_{k}<\delta^{2}$, then
	$\mathcal{H}_{1,k}$ is the subspace generated by
	$(e_{k},-r_{1,k}e_{k})$, where $\mathcal{A}$ acts as
	multiplication by $r_{1,k}$, and $\mathcal{H}_{2,k}$ is
	the subspace generated by $(e_{k},-r_{2,k}e_{k})$, where
	$\mathcal{A}$ acts as multiplication by $r_{2,k}$.
\end{itemize}

If we group together all subspaces corresponding to
eigenvalues of $\mathcal{A}$ with the same real part, we end
up with a decomposition of $\mathcal{H}$ of the form
\begin{equation}
	\mathcal{H}=\bigoplus_{r\in\mathcal{D}}\mathcal{H}_{r}.
	\label{decomp:hr}
\end{equation}

More precisely, $\calhr$ is the eigenspace of $\mathcal{A}$ 
relative to the real eigenvalue $r$ if $r\neq\delta$, while 
$\mathcal{H}_{\delta}$ is the closure of the space generated 
by all pairs of the form $(e_{k},0)$ and $(0,e_{k})$, with 
$k$ ranging over all indices for which the real part of 
$\lambda_{k}$ is greater than or equal to $\delta$.

We point out that (\ref{decomp:hr}) is a \emph{finite direct
sum of closed subspaces}.  In general it is not an orthogonal
sum, the only reason being that the pair of spaces originating
from each $\lambda_{k}<\delta^{2}$ are not orthogonal.

In this setting, a solution $u(t)$ to (\ref{eqn:lin-h}) is
\begin{itemize}
	\item $r_{0}$-pure if $(u(t),u'(t))\in
	\mathcal{H}_{r_{0}}$ for every $t\geq 0$,

	\item $r_{0}$-fast if $(u(t),u'(t))\in
	\bigoplus_{r>r_{0}}\mathcal{H}_{r}$  for every $t\geq 0$,

	\item $r_{0}$-slow if $(u(t),u'(t))\in
	\bigoplus_{r\leq r_{0}}\mathcal{H}_{r}$  for every $t\geq 0$.

\end{itemize}

\section{Statement of main results}\label{sec:statements}

Our first result provides a classification of all possible decay 
rates for solutions to the differential inequality (\ref{diff-ineq}).

\begin{thm}[Classification of decay rates]\label{thm:main-alternative}

	Let $H$ be a separable Hilbert space, let $A$ be a self-adjoint
	nonnegative operator on $H$ with dense domain $D(A)$, and let
	$u(t)$ be a global solution to the differential inequality
	(\ref{diff-ineq}) in the sense of Definition~\ref{defn:sol}.
	
	Let us assume that
	\begin{itemize}
		\item the spectrum of $A$ is a finite or countable set of
		eigenvalues without finite accumulation points,
	
		\item $K_{0}\geq 0$, $p>0$, and $q>0$,
	
		\item  $u$ is a decaying solution in the sense that
		\begin{equation}
			\lim_{t\to +\infty}|u(t)|_{D(A^{1/2})}=0.
			\label{hp:u-limit}
		\end{equation}
	\end{itemize}
		
	Then one and only one of the following statements apply.
	\begin{enumerate}
		\renewcommand{\labelenumi}{(\arabic{enumi})}
		
		\item \emph{(Null solution)} The solution is the zero-solution
		$u(t)\equiv 0$ for every $t\geq 0$.
		
		\item \emph{(Slow solutions)} There exist $T_{0}\geq 0$, and
		positive constants $M_{1}$ and $M_{2}$, such that
		\begin{equation}
			|u(t)|\geq \frac{M_{1}}{(1+t)^{1/p}}
			\quad\quad
			\forall t\geq T_{0},
			\label{th:slow}
		\end{equation}
		\begin{equation}
			|u'(t)|+|A^{1/2}u(t)|\leq M_{2}|u(t)|^{1+p}
			\quad\quad
			\forall t\geq T_{0}.
			\label{th:range}
		\end{equation}
		
		\item \emph{(Fast solutions)} There exist $r_{0}\in
		\mathcal{D}$, with $r_{0}>0$, and a nontrivial $r_{0}$-pure
		solution $v_{0}(t)$ to the linear homogeneous equation
		(\ref{eqn:lin-h}), such that
		\begin{equation}
			\lim_{t\to+\infty}\left(\left|u'(t)-v_{0}'(t)\right|+
			\left|u(t)-v_{0}(t)\right|_{D(A^{1/2})}\right)
			e^{\gamma t}=0
			\label{th:fast+}
		\end{equation}
		for every 
		\begin{equation}
			\gamma<\min\left\{\strut\beta_{0},(1+p)r_{0},
			(1+q)r_{0}\right\},
			\label{defn:delta}
		\end{equation}
		where
		\begin{equation}
			\beta_{0}:=\left\{
			\begin{array}{ll}
				\min\{r\in\mathcal{D}:r>r_{0}\} & 
				\mbox{if }r_{0}<\max\mathcal{D},  \\
				\noalign{\vspace{0.5ex}}
				+\infty & \mbox{otherwise}.
			\end{array}
			\right.
			\label{defn:beta-0}
		\end{equation}
	\end{enumerate}
\end{thm}

\begin{rmk}
	\begin{em}
		When the kernel of $A$ is non-trivial, a differential
		inequality such as~(\ref{diff-ineq}) does not guarantee that
		all its solutions in a neighborhood of the origin tend to 0.
		As an example, we can think to the ordinary differential
		equation $u''+u'=u^{3}+3u^{5}$, which has
		$u(t)=(1-2t)^{-1/2}$, and all its time-translations, among its
		solutions.  This is the reason why we need
		assumption~(\ref{hp:u-limit}).
		
		In other words, there might be coexistence of solutions that
		decay to 0 and solutions that do not decay, or even do not
		globally exist.  When this is the case, our result classifies
		all possible decay rates of decaying solutions, regardless of
		non-decaying ones.
	\end{em}
\end{rmk}

\begin{rmk}
	\begin{em}
		Concerning the null solution,
		Theorem~\ref{thm:main-alternative} implies that every
		solution which decays faster than $e^{-2\delta t}$ is
		actually the null solution.  This follows from the fact that
		all elements of $\mathcal{D}$ are less than or equal to
		$2\delta$.  In sharp contrast with the first order case
		(see~\cite{ggh:casc-par}), now there is a maximal possible
		decay rate for non-zero solutions.
		
	\end{em}
\end{rmk}

\begin{rmk}
	\begin{em}
		Concerning slow solutions, first of all we remark that in
		general (\ref{th:slow}) and (\ref{th:range}) involve a time
		$T_{0}\geq 0$.  This is because $u(t)$ could even vanish many
		times, of course not together with $u'(t)$, before becoming 
		slow ``eventually''. This point contrasts with the 
		parabolic case, where the corresponding estimates hold true 
		for every $t\geq 0$.
		
		For the rest, the slow regime is completely analogous to the
		first order case.  For example, estimates (\ref{th:slow}) and
		(\ref{th:range}) involve uniquely the exponent $p$ of the
		differential inequality (\ref{diff-ineq}), while $q$ is
		irrelevant provided it is positive.  Roughly speaking, this
		happens because slow solutions move closer and closer to the
		kernel of $A$, as suggested by the otherwise unnatural
		estimate (\ref{th:range}) in which $|A^{1/2}u(t)|$ is
		controlled by a higher power of $|u(t)|$.  Close to the kernel
		of $A$, the term $|A^{1/2}u(t)|$ can be neglected, and this
		justifies the disappearance of $q$ in the final decay rate.
		
		For the same reason, the slowness of $u(t)$ is due uniquely to
		its component with respect to $\ker(A)$.  Indeed, let us write
		$u(t)$ as the sum of its projection $P_{K}u(t)$ into
		$\ker(A)$, and its ``range component'' $u(t)-P_{K}u(t)$
		orthogonal to $\ker(A)$.  Since the operator $A$ is coercive
		when restricted to the range of $A$, estimate (\ref{th:range})
		implies that there exists a constant $c$ such that
		$$\left|u(t)-P_{K}u(t)\right|\leq
		c|A^{1/2}u(t)|\leq cM_{2}|u(t)|^{1+p}.$$
		
		Therefore, when $u(t)$ decays to 0, its range component always
		decays faster.  This is consistent with previous results
		(see~\cite{h:nodea}), and shows also that slow solutions
		cannot exist when the operator $A$ is coercive.
		
		Finally, we show that the exponent $(1+p)$ in (\ref{th:range})
		is optimal, both for $|u'(t)|$ and for $|A^{1/2}u(t)|$.
		Indeed, let us consider the case where $H=\re^{2}$, $p=2$, and
		the evolution problem reduces to the following system of
		ordinary differential equations 
		$$\left\{
		\begin{array}{l}
			x''(t)+x'(t)= -x^{3}(t)+3x^{5}(t),  \\
			\noalign{\vspace{0.5ex}}
			y''(t)+y'(t)+y(t)=x^{3}(t)-3x^{5}(t)+15x^{7}(t).
		\end{array}
		\right.$$
		
		A solution of this system is $x(t)=(1+2t)^{-1/2}$ and
		$y(t)=(1+2t)^{-3/2}$.  Therefore, in this case it turns out
		that $u(t)=(x(t),y(t))$ decays as $x(t)$, hence as $t^{-1/2}$,
		while both $u'(t)=(x'(t),y'(t))$ and $A^{1/2}u(t)=(0, y(t))$,
		decay as $t^{-3/2}$, hence as $|u(t)|^{1+p}$.

	\end{em}
\end{rmk}

\begin{rmk}
	\begin{em}
		Concerning fast solutions, we show that (\ref{defn:delta}) is
		optimal.  This can be seen by considering the case where
		$H=\re^{2}$ and the evolution problem reduces to the following
		system of ordinary differential equations 
		$$\left\{
		\begin{array}{l}
			x''(t)+x'(t)+(r_{0}-r_{0}^{2}) x(t)=0,  \\
			\noalign{\vspace{1ex}}
			y''(t)+y'(t)+(\beta_{0}-\beta_{0}^{2}) y(t)=
			|x(t)|^{1+p}+|x(t)|^{1+q},
		\end{array}
		\right.$$
		with parameters $0<r_{0}<\beta_{0}<1/2$.  It is not difficult
		to see that in this case
		$\mathcal{D}=\{r_{0},\beta_{0},1-\beta_{0},1-r_{0}\}$, and a
		solution of the first equation is $x(t)=e^{-r_{0}t}$.  At this
		point, solutions of the second equation can decay as $e^{-\eta
		t}$, where $\eta$ is the right-hand side
		of~(\ref{defn:delta}), or even as $te^{-\eta t}$ in case of
		resonance.
	\end{em}
\end{rmk}

The slow-fast alternative alone does not guarantee the existence of
both slow and fast solutions, and actually it does not guarantee the
existence of solutions at all.

As mentioned in the introduction, the existence of slow solutions was
addressed in~\cite{ggh:sol-lentes}, where an affirmative answer was
given assuming that the kernel of $A$ is nontrivial and that
$f(u)=-\nabla F(u)$ for a suitable nonnegative potential $F(u)$ such
that $\langle\nabla F(u),u\rangle\geq 0$ in a neighborhood of the
origin.  We refer for the details to Theorem~2.3
in~\cite{ggh:sol-lentes}.  Existence of slow solutions without this
gradient structure remains an open problem, as well as the slow-fast 
alternative when assumption (\ref{hp:u-limit}) is weakened by 
considering just the norm in $H$ instead of $D(A^{1/2})$.

In the next result we show the existence of families of fast solution 
with prescribed asymptotic profile. 

\begin{thm}[Existence of fast solutions]
	\label{thm:main-exponential}
	
	Let $H$ be a separable Hilbert space, and let $A$ be a self-adjoint
	nonnegative operator on $H$ with dense domain $D(A)$.  Let
	$f:B_{R_{0}}\to H$ be a function, with $R_{0}>0$ and $B_{R_{0}}$
	defined by (\ref{defn:ball}).
	
	Let us assume that
	\begin{enumerate}
		\renewcommand{\labelenumi}{(\roman{enumi})}
		
		\item the spectrum of $A$ is a finite or countable set of
		eigenvalues without finite accumulation points,
		
		\item there exist $p>0$ and $L\geq 0$ such that
		\begin{equation}
			|f(u)-f(v)|\leq L\left(
			|u|^{p}_{D(A^{1/2})}+|v|^{p}_{D(A^{1/2})}\right)|u-v|_{D(A^{1/2})}
			\label{hp:f-plip}
		\end{equation}
		for every $u$ and $v$ in $B_{R_{0}}$, and in addition
		\begin{equation}
			f(0)=0.
			\label{hp:f0}
		\end{equation}
	\end{enumerate}
	
	Then for every $r_{0}\in\mathcal{D}$, with $r_{0}>0$, there exists
	$\ep_{0}>0$ with the following property.  For every $r_{0}$-pure
	pair $(v_{0},v_{1})\in D(A^{1/2})\times H$, and for every
	$r_{0}$-fast pair $(z_{0},z_{1})\in D(A^{1/2})\times H$ such that
	\begin{equation}
		|v_{1}|+|v_{0}|_{D(A^{1/2})}+
		|z_{1}|+|z_{0}|_{D(A^{1/2})}\leq\ep_{0},
		\label{hp:fast}
	\end{equation}
	there exists an $r_{0}$-slow pair $(w_{0},w_{1})\in
	D(A^{1/2})\times H$ such that equation (\ref{eqn:semilin}) with
	initial data
	\begin{equation}
		u_{0}:=v_{0}+z_{0}+w_{0},
		\hspace{3em}
		u_{1}:=v_{1}+z_{1}+w_{1},
		\label{th:u-data}
	\end{equation}
	has a global solution $u(t)$ satisfying
	\begin{equation}
		\lim_{t\to +\infty}
		\left(|u'(t)-v'(t)|+
		|u(t)-v(t)|_{D(A^{1/2})}\right)e^{r_{0}t}=0,
		\label{th:limit}
	\end{equation}
	where $v(t)$ denotes the $r_{0}$-pure solution to the linear 
	homogeneous equation (\ref{eqn:lin-h}) with initial data 
	$(v_{0},v_{1})$.

\end{thm}

\begin{rmk}
	\begin{em}
		Once we know that $u(t)$ satisfies (\ref{th:limit}), then we 
		can always apply Theorem~\ref{thm:main-alternative} to it, 
		and improve the exponent $r_{0}$ to a better exponent.
	\end{em}
\end{rmk}

\begin{rmk}\label{rmk:exp-linear}
	\begin{em}
		Theorem~\ref{thm:main-exponential} can be seen as an existence
		result with a mix of conditions at $t=0$ and at $t=+\infty$.
		Indeed we fixed $(z_{0},z_{1})$, which at the end are the
		$r_{0}$-fast components of the initial condition
		$(u_{0},u_{1})$ (these components are not modified when adding
		$(v_{0},v_{1})$ and $(w_{0},w_{1})$), and through
		$(v_{0},v_{1})$ we also fixed $v(t)$, which at the end is the
		asymptotic profile of the solution as $t\to +\infty$.
		
		In the case of the linear homogeneous equation
		(\ref{eqn:lin-h}) it is not difficult to see that, for every
		$r_{0}$-pure asymptotic profile and every $r_{0}$-fast
		component of initial condition, there exists a (unique)
		solution satisfying these conditions.
		Theorem~\ref{thm:main-exponential} shows that the existence
		part is still true in the nonlinear case, at least if we limit
		ourselves to a neighborhood of the origin.
		
	\end{em}
\end{rmk}

\setcounter{equation}{0}
\section{Proofs}\label{sec:proofs}

\subsection{Estimates for non-homogeneous linear equations}\label{sec:lin}

The following result is the linear core of this paper.  In a few
words, we prove existence and uniqueness of a solution to the linear
non-homogeneous equation (\ref{eqn:lin-nh}) satisfying two conflicting
constraints.  The first one is that this solution decays almost as the
forcing term, the second one is that its initial data are as slow as
possible.

\begin{prop}[Special solution to the non-homogeneous equation]\label{prop:main-nh}
	Let $H$ be a separable Hilbert space, and let $A$ be a
	self-adjoint linear operator on $H$ with dense domain $D(A)$,
	whose spectrum $\sigma(A)$ is a finite or countable set of
	eigenvalues without finite accumulation points.  Let $g\in
	C^{0}([0,+\infty),H)$, and let us assume that there exist real
	numbers $K_{g}\geq 0$ and $\gamma_{0}>0$, with
	$\gamma_{0}\not\in\mathcal{D}$, such that
	\begin{equation}
		|g(t)|\leq K_{g}e^{-\gamma_{0}t}
		\quad\quad
		\forall t\geq 0.
		\label{hp:lin-nh}
	\end{equation}
	
	Then there exists a constant $\Gamma_{0}$, depending only on $\gamma_{0}$, 
	$\delta$ and $\sigma(A)$, for which the following statements 
	hold true.
	\begin{enumerate}
		\renewcommand{\labelenumi}{(\arabic{enumi})}
		\item  If $\gamma_{0}>\min\mathcal{D}$, and we set
		\begin{equation}
			\alpha_{0}:=\max\{r\in\mathcal{D}:r<\gamma_{0}\},
			\label{defn:alpha-0}
		\end{equation}
		then the non-homogeneous equation (\ref{eqn:lin-nh}) admits a
		unique solution $w(t)$ such that its initial condition is an
		$\alpha_{0}$-slow pair and
		\begin{equation}
			\lim_{t\to+\infty}
			\left(|w'(t)|+|w(t)|_{D(A^{1/2})}\right)e^{\alpha_{0}t}=0.
			\label{th:lin-nh-lim}
		\end{equation}
		
		Moreover, this solution satisfies the stronger decay estimate
		\begin{equation}
			|w'(t)|+|w(t)|_{D(A^{1/2})}\leq 
			\Gamma_{0}K_{g}e^{-\gamma_{0}t}
			\quad\quad
			\forall t\geq 0.
			\label{th:lin-nh}
		\end{equation}

		\item  If $\gamma_{0}<\min\mathcal{D}$, then the solution to 
		the non-homogeneous equation (\ref{eqn:lin-nh}) with initial 
		data $w(0)=w'(0)=0$ satisfies (\ref{th:lin-nh}).
	
	\end{enumerate}
\end{prop}

\paragraph{\textmd{\textit{Proof}}}

\subparagraph{\textmd{\textit{Uniqueness}}}

In the case
$\gamma_{0}<\min\mathcal{D}$ uniqueness is trivial because both
initial data are given.  In the case $\gamma_{0}>\min\mathcal{D}$, let
$w_{1}(t)$ and $w_{2}(t)$ be two solutions, and let
$v(t):=w_{1}(t)-w_{2}(t)$ denote their difference.  Clearly $v(t)$ is
a solution to the corresponding homogeneous equation.  Moreover, it is
$\alpha_{0}$-slow because its initial data are the difference of two
$\alpha_{0}$-slow pairs, and it satisfies
\begin{equation}
	\lim_{t\to+\infty}
	\left(|v'(t)|+|v(t)|_{D(A^{1/2})}\right)e^{\alpha_{0}t}=0
	\label{est:diff-lim}
\end{equation}
because the same is true for $w_{1}(t)$ and $w_{2}(t)$. On the other 
hand, the unique $\alpha_{0}$-slow solution for which 
(\ref{est:diff-lim}) holds true is the null solution, as pointed out 
in Remark~\ref{rmk:lin-h-est}. This proves 
that $w_{1}(t)=w_{2}(t)$.

\subparagraph{\textmd{\textit{Estimates in the product space}}}

Let us interpret the non-homogeneous equation (\ref{eqn:lin-nh}) as a
first order system, as we did in section~\ref{sec:defn-system} for the
homogeneous equation.  Setting $W(t):=(w(t),w'(t))$ it turns out that
$w(t)$ solves (\ref{eqn:lin-nh}) in $H$ if and only if $W(t)$ solves
the first order system
\begin{equation}
	W'(t)+\mathcal{A}W(t)=(0,g(t))
	\label{system-W}
\end{equation}
in the product space $\mathcal{H}:=D(A^{1/2})\times H$, with the 
operator $\mathcal{A}$ defined by (\ref{defn:cal-A}).

Let us assume for simplicity that
\begin{equation}
	\min\mathcal{D}<\gamma_{0}<\max\mathcal{D},
	\label{hp:gamma-0}
\end{equation}
and let us define $\alpha_{0}$ as in (\ref{defn:alpha-0}), and 
$\beta_{0}$ as the smallest element of $\mathcal{D}$ greater than 
$\gamma_{0}$. Let us write
\begin{equation}
	\calh=\mathcal{H}_{-}\oplus\mathcal{H}_{+},
	\label{decomp:H+-}
\end{equation}
where 
$$\mathcal{H}_{-}:=
\bigoplus_{\stackrel{\scriptstyle{r\in\mathcal{D}\rule[-0.5ex]{0pt}{0pt}}}{r<\gamma_{0}}}
\mathcal{H}_{r}=
\bigoplus_{\stackrel{\scriptstyle{r\in\mathcal{D}\rule[-0.5ex]{0pt}{0pt}}}{r\leq\alpha_{0}}}
\mathcal{H}_{r}, \hspace{4em} \mathcal{H}_{+}:=
\bigoplus_{\stackrel{\scriptstyle{r\in\mathcal{D}\rule[-0.5ex]{0pt}{0pt}}}{r>\gamma_{0}}}
\mathcal{H}_{r}=
\bigoplus_{\stackrel{\scriptstyle{r\in\mathcal{D}\rule[-0.5ex]{0pt}{0pt}}}{r\geq\beta_{0}}}
\mathcal{H}_{r}.$$

Since (\ref{decomp:hr}) is a direct sum of closed subspaces, the
projection onto each $\calh_{r}$ is continuous with a norm depending 
only on $\mathcal{D}$.  As a consequence, the
projections $P_{-}$ and $P_{+}$ onto $\mathcal{H}_{-}$ and
$\mathcal{H}_{+}$, respectively, are continuous, hence there exist two
constants $\Gamma_{1}$ and $\Gamma_{2}$, depending only on
$\mathcal{D}$, such that
\begin{equation}
	\left|P_{+}(x,y)\right|_{\calh}\leq
	\Gamma_{1}\left|(x,y)\right|_{\calh}
	\quad\quad
	\forall (x,y)\in\calh,
	\label{est:p+}
\end{equation}
\begin{equation}
	\left|P_{-}(x,y)\right|_{\calh}\leq
	\Gamma_{2}\left|(x,y)\right|_{\calh}
	\quad\quad
	\forall (x,y)\in\calh,
	\label{est:p-}
\end{equation}
where of course $|(x,y)|_{\calh}^{2}:=|x|^{2}_{D(A^{1/2})}+|y|^{2}$.

Now let us consider the semigroup $S(t)$ generated on $\calh$ by the 
first order system (\ref{system}). From the explicit solutions (\ref{sol:uk-sup}) through 
(\ref{sol:uk-sub}) it follows that
\begin{equation}
	\left|S(t)(x,y)\right|_{\calh}\leq
	\Gamma_{3}(1+t)e^{-\beta_{0}t}\left|(x,y)\right|_{\calh}
	\quad\quad
	\forall t\geq 0,\quad\forall(x,y)\in\mathcal{H}_{+},
	\label{est:S+}
\end{equation}
where $\Gamma_{3}$ depends only on $\mathcal{D}$ and $\sigma(A)$.

The operators $S(t)$ of the semigroup are invertible, and $S(-t)$ 
corresponds to solving the second order equation
$$u''(t)-2\delta u'(t)+Au(t)=0,$$
in which we just reversed the sign of the damping term.  Explicit solutions
in this case are analogous to (\ref{sol:uk-sup}) through
(\ref{sol:uk-sub}), just with the opposite sign in the argument of all
exponentials. As a consequence, $S(-t)$ can be estimates in 
$\mathcal{H}_{-}$ as follows:
\begin{equation}
	\left|S(-t)(x,y)\right|_{\calh}\leq
	\Gamma_{4}(1+t)e^{\alpha_{0}t}\left|(x,y)\right|_{\calh}
	\quad\quad
	\forall t\geq 0,\quad\forall(x,y)\in\mathcal{H}_{-},
	\label{est:S-}
\end{equation}
where $\Gamma_{4}$ depends only on $\mathcal{D}$ and $\sigma(A)$.

When assumption (\ref{hp:gamma-0}) is not satisfied,
the situation is even simpler.  If $\gamma_{0}>\max\mathcal{D}$, then
$\mathcal{H}_{-}=\calh$ and $\mathcal{H}_{+}=\{0\}$, and 
(\ref{est:S-}) still holds true with $\alpha_{0}$ given by
(\ref{defn:alpha-0}). If $\gamma_{0}<\min\mathcal{D}$,
then $\mathcal{H}_{-}=\{0\}$ and $\mathcal{H}_{+}=\calh$, and 
(\ref{est:S+}) still holds true with the same $\beta_{0}$.

\subparagraph{\textmd{\textit{Existence}}}

We are now ready to prove existence of the solution with the required 
properties. Let us set
\begin{equation}
	W_{+}(t):=\int_{0}^{t}S(t-s)P_{+}(0,g(s))\,ds,
	\label{defn:w+}
\end{equation}
\begin{equation}
	W_{-}(t):=-\int_{t}^{+\infty}S(t-s)P_{-}(0,g(s))\,ds.
	\label{defn:w-}
\end{equation}

Let us assume for a while that $W_{\pm}(t)$ are well-defined, namely 
that the integrals are convergent. Then they are solutions to
$$W_{\pm}'(t)+\mathcal{A}W_{\pm}(t)=P_{\pm}(0,g(t)),$$
hence their sum $W(t):=W_{+}(t)+W_{-}(t)$ is a solution to
(\ref{system-W}), which is equivalent to saying that $W(t)$ is of the
form $(w(t),w'(t))$ for some solution $w(t)$ to (\ref{eqn:lin-nh}).
Moreover, the initial condition $W(0)$ is $\alpha_{0}$-slow, because
its $\alpha_{0}$-fast component $W_{+}(0)$ vanishes due to
(\ref{defn:w+}) with $t=0$.

Therefore, we are left to proving that $W_{\pm}(t)$ are well-defined 
and satisfy suitable decay estimates. 
Let us start with $W_{+}(t)$.  Since $t-s\geq 0$ in (\ref{defn:w+}),
from (\ref{est:S+}), $(\ref{est:p+})$, and (\ref{hp:lin-nh}) it
follows that
\begin{eqnarray*}
	\left|S(t-s)P_{+}(0,g(s))\right| & \leq & 
	\Gamma_{3}(1+t-s)e^{-\beta_{0}(t-s)}\cdot
	\left|P_{+}(0,g(s))\right|_{\calh}\\
	 & \leq & \Gamma_{3}(1+t-s)e^{-\beta_{0}(t-s)}\cdot
	 \Gamma_{1} |g(s)|\\
	 & \leq & \Gamma_{3}(1+t-s)e^{-\beta_{0}(t-s)}\cdot
	 \Gamma_{1} K_{g}e^{-\gamma_{0}s}  \\
	 & = & \Gamma_{1}\Gamma_{3}K_{g}\cdot
	 e^{-\gamma_{0}t}\cdot (1+t-s)e^{-(\beta_{0}-\gamma_{0})(t-s)},
\end{eqnarray*}
and hence
\begin{equation}
	|W_{+}(t)|\leq\Gamma_{1}\Gamma_{3}K_{g}\cdot
	e^{-\gamma_{0}t}\cdot 
	\int_{0}^{t}(1+t-s)e^{-(\beta_{0}-\gamma_{0})(t-s)}\,ds\leq
	\Gamma_{5}K_{g}e^{-\gamma_{0}t}.
	\label{est:w+}
\end{equation}

Now let us consider $W_{-}(t)$.  Since $t-s\leq 0$ in (\ref{defn:w-}),
from (\ref{est:S-}), $(\ref{est:p-})$, and (\ref{hp:lin-nh}) it
follows that
\begin{eqnarray*}
	\left|S(t-s)P_{-}(0,g(s))\right| & \leq & 
	\Gamma_{4}(1+s-t)e^{\alpha_{0}(s-t)}\cdot
	\left|P_{-}(0,g(s))\right|_{\calh}\\
	 & \leq & \Gamma_{4}(1+s-t)e^{\alpha_{0}(s-t)}\cdot
	 \Gamma_{2} |g(s)|\\
	 & \leq & \Gamma_{4}(1+s-t)e^{\alpha_{0}(s-t)}\cdot
	 \Gamma_{2} K_{g}e^{-\gamma_{0}s}  \\
	 & = & \Gamma_{2}\Gamma_{4}K_{g}\cdot
	 e^{-\gamma_{0}t}\cdot (1+s-t)e^{-(\gamma_{0}-\alpha_{0})(s-t)}.
\end{eqnarray*}

Since $\gamma_{0}>\alpha_{0}$, this proves that the integral in 
(\ref{defn:w-}) converges. In addition, it turns out that
\begin{equation}
	|W_{-}(t)|\leq\Gamma_{2}\Gamma_{4}K_{g}\cdot
	e^{-\gamma_{0}t}\cdot 
	\int_{t}^{+\infty}(1+s-t)e^{-(\gamma_{0}-\alpha_{0})(s-t)}\,ds\leq
	\Gamma_{6}K_{g}e^{-\gamma_{0}t}.
	\label{est:w-}
\end{equation}
Summing (\ref{est:w+}) and (\ref{est:w-}) it follows that
$$|W(t)|_{\calh}\leq(\Gamma_{5}+\Gamma_{6})K_{g}
e^{-\gamma_{0}t}
\quad\quad
\forall t\geq 0,$$
which is equivalent to (\ref{th:lin-nh}).\qed

\begin{rmk}
	\begin{em}
		Proposition~\ref{prop:main-nh} can also be proved without 
		relying on the product space $\calh$. It is enough to define 
		the components $w_{k}(t)$ of the solution as suitable 
		integrals involving the components $g_{k}(t)$ of the forcing 
		term and the fundamental solutions of the homogeneous equation.
		This requires to distinguish several cases. For example, when 
		$\lambda_{k}>\delta^{2}$ the component is given by
		$$w_{k}(t):=\frac{1}{(\lambda_{k}-\delta^{2})^{1/2}}
		\int_{I}e^{-\delta(t-s)}
		\sin\left((\lambda_{k}-\delta^{2})^{1/2}(t-s)\right)
		g_{k}(s)\,ds,$$
		where the integration region is $I:=[0,t]$ if 
		$\delta>\gamma_{0}$, and $I:=[t,+\infty)$ if 
		$\delta<\gamma_{0}$.
	\end{em}
\end{rmk}

\subsection{Energies and Dirichlet quotients}\label{sec:energies}

In this section we introduce the energies which we are going to
exploit in the proof of Theorem~\ref{thm:main-alternative}.  To begin
with, we consider the usual ``hyperbolic energy'' $E(t)$ defined in
(\ref{defn:E}), and the following \emph{generalized Dirichlet
quotient}
\begin{equation}
	G_{d}(t):=\frac{|u'(t)|^{2}+|A^{1/2}u(t)|^{2}}{|u(t)|^{2+d}}
	=\frac{E(t)}{|u(t)|^{2+d}},
	\label{defn:Qd}
\end{equation}
defined for every $d\geq 0$ provided that $u(t)\neq 0$.  This is the
hyperbolic version of the quotient
$$Q_{d}(t):=\frac{|A^{1/2}u(t)|^{2}}{|u(t)|^{2+d}},$$
introduced in~\cite{ggh:casc-par} when dealing with the semilinear
parabolic problem.

We need also a modified version of the hyperbolic energy
and of the generalized Dirichlet quotient, which we define in the
following way.  Let $Q:H\to H$ denote the orthogonal projection onto
$\ker(A)^{\perp}$.  Thanks to our assumptions on the spectrum of $A$,
there exists a constant $\nu>0$ such that
\begin{equation}
	|Qu|^{2}\leq\frac{1}{\nu}|A^{1/2}u|^{2}
	\quad\quad
	\forall u\in D(A^{1/2}).
	\label{est:Qu}
\end{equation}

More precisely, we can take $\nu$ equal to any positive number if $A$ 
is the null operator, and $\nu$ equal to the smallest positive 
eigenvalue of $A$ otherwise. Now let us set
\begin{equation}
	\mu:=\min\left\{\frac{1}{2},\frac{\nu}{2},
	\frac{\delta}{2},\frac{\nu}{5\delta}\right\},
	\label{defn:gamma}
\end{equation}
and let us define the \emph{modified hyperbolic energy}
\begin{equation}
	\Ehat(t):=|u'(t)|^{2}+|A^{1/2}u(t)|^{2}+
	2\mu\langle u'(t),Qu(t)\rangle,
	\label{defn:Fhat}
\end{equation}
and finally the \emph{modified generalized Dirichlet quotient}
\begin{equation}
	\Ghat_{d}(t):=\frac{|u'(t)|^{2}+|A^{1/2}u(t)|^{2}+
	2\mu\langle u'(t),Qu(t)\rangle}{|u(t)|^{2+d}}=
	\frac{\Ehat(t)}{|u(t)|^{2+d}}.
	\label{defn:Ghatd}
\end{equation}

Next result provides estimates on these quantities and their
time-derivatives in the case where $u(t)$ is a solution of a
non-homogeneous linear equation such as~(\ref{eqn:lin-nh}).

\begin{lemma}[Energies and generalized Dirichlet quotients]
	Let $H$ be a separable Hilbert space, and let $A$ be a self-adjoint
	nonnegative operator on $H$ with dense domain $D(A)$.  Let $g\in
	C^{0}([0,+\infty),H)$, and let $u(t)$ be a solution
	to~(\ref{eqn:lin-nh}) in $[0,+\infty)$ in the sense of
	Theorem~\ref{thmbibl:linear}.
	
	Then the energies $E(t)$, $\Ehat(t)$, $G_{d}(t)$, $\Ghat_{d}(t)$
	satisfy the following estimates.
	\begin{enumerate}
		\renewcommand{\labelenumi}{(\arabic{enumi})}
		\item  It turns out that
		\begin{equation}
			\frac{1}{2}E(t)\leq \Ehat(t)\leq 2E(t)
			\quad\quad
			\forall t\geq 0.
			\label{th:Fhat-F}
		\end{equation}
		
		Moreover, the function $\Ehat(t)$ is of class $C^{1}$ in 
		$[0,+\infty)$ and
		\begin{equation}
			\Ehat'(t)\leq -\frac{\mu}{2}\Ehat(t)+
			\frac{2}{\delta}|g(t)|^{2}
			\quad\quad
			\forall t\geq 0.
			\label{th:deriv-Fhat}
		\end{equation}
	
		\item Let us assume that $u(t)\neq 0$ for every $t$ in some 
		time-interval $(a,b)$.
		Then the generalized Dirichlet quotients are well defined in
		$(a,b)$, and for every $d\geq 0$ it turns out that
		\begin{equation}
			\frac{1}{2}G_{d}(t)\leq
			\Ghat_{d}(t)\leq
			2G_{d}(t)
			\quad\quad
			\forall t\in(a,b).
			\label{th:Ghat-G}
		\end{equation}
		
		Moreover, the function $\Ghat_{d}(t)$ is of class $C^{1}$ in 
		$(a,b)$ and
		\begin{equation}
			\Ghat_{d}'(t)\leq -\frac{\mu}{2}\Ghat_{d}(t)+
			\frac{2}{\delta}\frac{|g(t)|^{2}}{|u(t)|^{2+d}}+			
			(2+d)|u(t)|^{d/2}\cdot
			\left[G_{d}(t)\right]^{1/2}\cdot\Ghat_{d}(t)
			\label{th:deriv-Ghat}
		\end{equation}
		for every $t\in(a,b)$.
	\end{enumerate}
\end{lemma}

\paragraph{\textmd{\textit{Proof}}}

From (\ref{est:Qu}) we obtain that
$$\left|2\langle u'(t),Qu(t)\rangle\right|\leq
2|u'(t)|\cdot|Qu(t)|\leq
|u'(t)|^{2}+|Qu(t)|^{2}\leq
|u'(t)|^{2}+\frac{1}{\nu}|A^{1/2}u(t)|^{2}.$$

Since $\mu\leq 1/2$ and $\mu\leq\nu/2$, it follows that
$$\left|2\mu\langle u'(t),Qu(t)\rangle\right|\leq
\frac{1}{2}|u'(t)|^{2}+\frac{1}{2}|A^{1/2}u(t)|^{2},$$
which proves both (\ref{th:Fhat-F}) and (\ref{th:Ghat-G}).

The time-derivative of $\Ehat(t)$ is 
\begin{eqnarray}
	\Ehat'(t) & = & -4\delta|u'(t)|^{2}-
	2\mu\langle A^{1/2}u(t),A^{1/2}Qu(t)\rangle+
	2\mu\langle u'(t),Qu'(t)\rangle 
	\nonumber   \\
	\noalign{\vspace{1ex}}
	 & & \mbox{}-4\mu\delta\langle u'(t),Qu(t)\rangle +
	 2\langle u'(t),g(t)\rangle+
	 2\mu\langle Qu(t),g(t)\rangle.
	 \label{deriv:Fhat}
\end{eqnarray}

Let $I_{1}$, \ldots, $I_{6}$ denote the six terms in the right-hand 
side, which we now estimate separately. From the definition of $Q$ it 
follows that
\begin{equation}
	I_{2}=-2\mu|A^{1/2}u(t)|^{2}
	\label{est:Fhat-I2}
\end{equation}
and $\langle u'(t),Qu'(t)\rangle\leq|u'(t)|^{2}$. Since 
$2\mu\leq\delta$, we deduce that
\begin{equation}
	I_{3}\leq\delta|u'(t)|^{2}.
	\label{est:Fhat-I3}
\end{equation}

Since
$$I_{4}\leq 4\mu\delta|u'(t)|\cdot|Qu(t)|\leq
2\mu\delta\left(
\frac{|u'(t)|^{2}}{2\mu}+2\mu|Qu(t)|^{2}\right),$$
from (\ref{est:Qu}) it follows that
\begin{equation}
	I_{4}\leq\delta|u'(t)|^{2}+
	\frac{4\mu^{2}\delta}{\nu}|A^{1/2}u(t)|^{2}.
	\label{est:Fhat-I4}
\end{equation}

As for the last two terms, it turns out that
\begin{equation}
	I_{5}\leq\delta|u'(t)|^{2}+\frac{1}{\delta}|g(t)|^{2},
	\label{est:Fhat-I5}
\end{equation}
and
$$I_{6}\leq 2\mu|Qu(t)|\cdot|g(t)|\leq
\mu\left(\mu\delta|Qu(t)|^{2}+
\frac{1}{\mu\delta}|g(t)|^{2}\right),$$
so that from (\ref{est:Qu}) it follows that
\begin{equation}
	I_{6}\leq\frac{\mu^{2}\delta}{\nu}|A^{1/2}u(t)|^{2}+
	\frac{1}{\delta}|g(t)|^{2}.
	\label{est:Fhat-I6}
\end{equation}

Plugging (\ref{est:Fhat-I2}) through (\ref{est:Fhat-I6}) into
(\ref{deriv:Fhat}), we obtain that 
$$\Ehat'(t)\leq -\delta|u'(t)|^{2}
-\mu\left(2-\frac{5\mu\delta}{\nu}\right)|A^{1/2}u(t)|^{2}
+\frac{2}{\delta}|g(t)|^{2}.$$

Keeping into account that $\mu\leq\delta$ and 
$5\mu\delta\leq\nu$, we conclude that
$$\Ehat'(t)\leq-\mu
\left(|u'(t)|^{2}+|A^{1/2}u(t)|^{2}\right)
+\frac{2}{\delta}|g(t)|^{2}.$$

At this point (\ref{th:deriv-Fhat}) follows from (\ref{th:Fhat-F}).

It remains to compute the time-derivative of $\Ghat_{d}(t)$, which 
turns out to be
\begin{equation}
	\Ghat_{d}'(t)=\frac{\Ehat'(t)}{|u(t)|^{2+d}} -
	(2+d)\frac{\langle u'(t),u(t)\rangle}{|u(t)|^{2}}
	\cdot\Ghat_{d}(t).
	\label{deriv:Ghat}
\end{equation}

From (\ref{th:deriv-Fhat}) it follows that
\begin{equation}
	\frac{\Ehat'(t)}{|u(t)|^{2+d}}\leq
	-\frac{\mu}{2}\frac{\Ehat(t)}{|u(t)|^{2+d}}+
	\frac{2}{\delta}\frac{|g(t)|^{2}}{|u(t)|^{2+d}},
	\label{est:deriv-Ghat-1}
\end{equation}
and from the definition of $G_{d}(t)$ it follows that
\begin{equation}
	\frac{\langle u'(t),u(t)\rangle}{|u(t)|^{2}}\leq
	\frac{|u'(t)|}{|u(t)|^{1+d/2}}\cdot|u(t)|^{d/2}\leq
	\left[G_{d}(t)\right]^{1/2}|u(t)|^{d/2}.
	\label{est:deriv-Ghat-2}
\end{equation}

Plugging (\ref{est:deriv-Ghat-1}) and (\ref{est:deriv-Ghat-2}) into 
(\ref{deriv:Ghat}), we obtain (\ref{th:deriv-Ghat}).\qed

\subsection{Proof of Theorem~\ref{thm:main-alternative}}

Let us describe the scheme and the heuristic ideas behind the proof
before entering into details.  In the first section of the proof we
get rid of the null solution.  Indeed we prove that
$(u(T),u'(T))=(0,0)$ for some $T\geq 0$ if and only if
$(u(t),u'(t))=(0,0)$ for every $t\geq 0$.  This is a result of forward
and backward uniqueness of the null solution.  After proving it, we
can assume that
\begin{equation}
	(u(t),u'(t))\neq(0,0)
	\quad\quad
	\forall t\geq 0.
	\label{hp:u-neq-0}
\end{equation}

In the second section of the proof we assume that there exist a
constant $c_{1}>0$ and a sequence $t_{n}\to +\infty$ such that
\begin{equation}
	|u'(t_{n})|^{2}+|A^{1/2}u(t_{n})|^{2}\leq c_{1}|u(t_{n})|^{2+p}
	\quad\quad
	\forall n\in\n.
	\label{hp:liminf}
\end{equation}

Under this assumption, we prove that a similar estimate holds true 
eventually, namely there exists $n_{0}\in\n$ such that
\begin{equation}
	|u'(t)|^{2}+|A^{1/2}u(t)|^{2}\leq 4c_{1}|u(t)|^{2+p}
	\quad\quad
	\forall t\geq t_{n_{0}}.
	\label{hp:liminf-glob}
\end{equation}

This is not yet (\ref{th:range}), but in any case it shows that
$|A^{1/2}u(t)|$ decays faster than $|u(t)|$.  As already pointed out,
this means that $u(t)$ moves closer and closer to the kernel of $A$,
and suggests that the terms with $Au(t)$ and $A^{1/2}u(t)$ in the
differential inequality (\ref{diff-ineq}) can be neglected.  Moreover,
since we expect solutions decaying as negative powers of $t$, it seems
reasonable to neglect also $u''(t)$, which for negative powers of $t$
decays faster than $u'(t)$.  With this ansatz, the second order
differential inequality (\ref{diff-ineq}) has become the first order
differential inequality $|u'(t)|\leq K_{0}|u(t)|^{1+p}$, whose nonzero
solutions are slow in the sense of (\ref{th:slow}).
The formal proof requires a sharp analysis of the Dirichlet quotients
of section~\ref{sec:energies}, first with $d:=p$ and then with
$d:=2p$.

In the third and last section of the proof we are left with the case where
(\ref{hp:liminf}) is false for every constant $c_{1}$ and every
sequence $t_{n}\to+\infty$.  This easily implies the existence of
$T_{1}\geq 0$ such that
\begin{equation}
	|u(t)|^{2+p}\leq |u'(t)|^{2}+|A^{1/2}u(t)|^{2}
	\quad\quad
	\forall t\geq T_{1}.
	\label{est:u<Au}
\end{equation}

In this case we are not allowed to ignore the operator $A$, but we can
neglect the right-hand side of (\ref{diff-ineq}) because the exponents
are larger than one.  Therefore, a good approximation of
(\ref{diff-ineq}) is now the linear homogeneous equation
(\ref{eqn:lin-h}), whose solutions decay exponentially with possible
rates corresponding to elements of $\mathcal{D}$.  The formal proof
requires several steps.  First of all, we provide exponential
estimates from below and from above with non-optimal rates.  Then we
identify the exact rate, and finally we prove that (\ref{th:fast+})
holds true. The basic tool in this part of the proof is 
Proposition~\ref{prop:main-nh}.

We point out that the exponent $2+p$ is non-optimal both in
(\ref{hp:liminf}) and in the opposite estimate (\ref{est:u<Au}).
Indeed, a posteriori it turns out that (up to multiplicative
constants) $|A^{1/2}u|\leq|u|^{1+p}$ in the case of slow solutions,
and $|A^{1/2}u|\sim|u|$ in the case of fast solutions, so that the
right exponents would be $2+2p$ in the slow regime and $2$ in the fast
regime.  Nevertheless, the intermediate exponent $2+p$ acts as a
threshold separating the two different regimes, and leaving enough
room on both sides to perform our estimates.

\subsubsection*{Non-trivial solutions never vanish in the phase space}

We prove that either $(u(t),u'(t))=(0,0)$ for every $t\geq 0$ or 
$(u(t),u'(t))\neq(0,0)$ for every $t\geq 0$.

To this end, we consider the energy
$$F(t):=|u'(t)|^{2}+|A^{1/2}u(t)|^{2}+|u(t)|^{2}.$$

Its time-derivative is
$$F'(t)=
-4\delta|u'(t)|^{2}+2\langle u'(t),g(t)\rangle+
2\langle u'(t),u(t)\rangle,$$
with $g(t)$ as in Definition~\ref{defn:sol}.  From assumption
(\ref{hp:main-g}) it follows that
\begin{equation}
	|g(t)|^{2}\leq 2K_{0}^{2}\left(
	|u(t)|^{2(1+p)}+|A^{1/2}u(t)|^{2(1+q)}\right)
	\quad\quad
	\forall t\geq 0,
	\label{hp:main-g-bis}
\end{equation}
and hence in this case
$$|g(t)|^{2}\leq 2K_{0}^{2}\left( [F(t)]^{1+p}+[F(t)]^{1+q}\right)
\quad\quad
\forall t\geq 0.$$

Thus it follows that
\begin{eqnarray*}
	|F'(t)| & \leq & 4\delta|u'(t)|^{2}+|u'(t)|^{2}+|g(t)|^{2}
	+|u'(t)|^{2}+|u(t)|^{2}\\
	\noalign{\vspace{0.5ex}}
	 & \leq & (4\delta +2)F(t)+
	 2K_{0}^{2}\left( [F(t)]^{1+p}+[F(t)]^{1+q}\right)
\end{eqnarray*}
for every $t\geq 0$.  The exponents of $F(t)$ in the right-hand side
are all greater than or equal to 1.  Therefore, this differential
inequality guarantees that either $F(t)=0$ for every $t\geq 0$ or
$F(t)>0$ for every $t\geq 0$, which is equivalent to what we had to
prove.

\subsubsection*{Slow solutions}

In this second part of the proof we consider the case where
(\ref{hp:u-neq-0}) holds true and $u(t)$ satisfies (\ref{hp:liminf})
for some $c_{1}> 0$ and some sequence $t_{n}\to +\infty$.

\paragraph{\textmd{\textit{Main estimate}}}

Let $\nu$ be the constant which appears in (\ref{est:Qu}), and let
$\mu$ be defined as in (\ref{defn:gamma}).  Due to assumption
(\ref{hp:u-limit}), there exists $n_{0}\in\n$ such that
\begin{equation}
	2(2+p)\sqrt{c_{1}}|u(t)|^{p/2}\leq\frac{\mu}{4}
	\quad\quad
	\forall t\geq t_{n_{0}},
	\label{hp:n0-1}
\end{equation}
\begin{equation}
	\frac{4K_{0}^{2}}{\delta}\left(
	|u(t)|^{p}+(4c_{1})^{1+q}|u(t)|^{(2+p)q}\right)\leq
	\frac{\mu}{4}c_{1}
	\quad\quad
	\forall t\geq t_{n_{0}}.
	\label{hp:n0-2}
\end{equation}

We claim that (\ref{hp:liminf-glob}) holds true with this choice of 
$n_{0}$, and that in addition
\begin{equation}
	|u(t)|>0
	\quad\quad
	\forall t\geq t_{n_{0}}.
	\label{th:u>0}
\end{equation}

To this end, let us consider the generalized Dirichlet
quotient~(\ref{defn:Qd}), and its modified version (\ref{defn:Ghatd}),
with $d:=p$.  To begin with, we observe that $u(t_{n_{0}})\neq 0$,
because if not we could deduce from (\ref{hp:liminf}) that
$(u(t_{n_{0}}),u'(t_{n_{0}}))=(0,0)$, and this would contradict
assumption (\ref{hp:u-neq-0}).  As a consequence, $G_{p}(t)$ and
$\Ghat_{p}(t)$ are defined at least in a neighborhood of $t_{n_{0}}$.
Moreover, since $G_{p}(t_{n_{0}})\leq c_{1}<4c_{1}$, by a continuity
argument it follows that $G_{p}(t)<4c_{1}$ in a suitable neighborhood
of $t_{n_{0}}$.  Let us set
$$S:=\sup\left\{t>t_{n_{0}}:
|u(\tau)|>0\mbox{ and }G_{p}(\tau)\leq 4c_{1}\quad
\forall\tau\in[t_{n_{0}},t]\right\},$$
so that (\ref{hp:liminf-glob}) and (\ref{th:u>0}) are now equivalent
to showing that $S=+\infty$.

Let us assume by contradiction that $S<+\infty$. By the maximality of 
$S$, this means that either $u(S)=0$ or $G_{p}(S)=4c_{1}$. Now we 
show that both options are impossible. In order to exclude the first 
one, we observe that
$$|u'(t)|^{2}+|A^{1/2}u(t)|^{2}\leq 4c_{1}|u(t)|^{2+p}
\quad\quad
\forall t\in[t_{n_{0}},S).$$

If $|u(S)|=0$, then letting $t\to S^{-}$ we deduce that also 
$|u'(S)|=0$, which contradicts again (\ref{hp:u-neq-0}).

It remains to exclude that $G_{p}(S)=4c_{1}$. Setting $d:=p$ in 
(\ref{th:deriv-Ghat}) we obtain that
$$\Ghat_{p}'(t)\leq
-\frac{\mu}{2}\Ghat_{p}(t)+
(2+p)|u(t)|^{p/2}\cdot[G_{p}(t)]^{1/2}\cdot\Ghat_{p}(t)+
\frac{2}{\delta}\frac{|g(t)|^{2}}{|u(t)|^{2+p}}.$$

Therefore, since (\ref{hp:main-g-bis}) implies that
$$|g(t)|^{2}\leq 2K_{0}^{2}|u(t)|^{2+p}\left(
|u(t)|^{p}+[G_{p}(t)]^{1+q}|u(t)|^{(2+p)q}\right),$$
we find that $\Ghat_{p}'(t)$ is less than or equal to
$$-\left(
\frac{\mu}{2}-(2+p)|u(t)|^{p/2}[G_{p}(t)]^{1/2}\right)\Ghat_{p}(t)+
\frac{4K_{0}^{2}}{\delta}\left(|u(t)|^{p}+
|u(t)|^{(2+p)q}\left[G_{p}(t)\right]^{1+q}\right).$$

If we keep into account that $G_{p}(t)\leq 4c_{1}$ for every 
$t\in[t_{n_{0}},S)$, and the smallness assumptions (\ref{hp:n0-1}) 
and (\ref{hp:n0-2}), we conclude that
$$\Ghat_{p}'(t)\leq -\frac{\mu}{4}\Ghat_{p}(t)+
\frac{\mu}{4}c_{1}
\quad\quad
\forall t\in[t_{n_{0}},S),$$
and hence
$$\Ghat_{p}(t)\leq
\left(\Ghat_{p}(t_{n_{0}})-c_{1}\right)
\exp\left(-\frac{\mu}{4}(t-t_{n_{0}})\right)+c_{1}
\quad\quad
\forall t\in[t_{n_{0}},S).$$

Now from (\ref{th:Ghat-G}) and (\ref{hp:liminf}) we know that 
$\Ghat_{p}(t_{n_{0}})\leq 2G_{p}(t_{n_{0}})\leq 2c_{1}$, so that
$$\Ghat_{p}(t)\leq c_{1}\left(1+
\exp\left(-\frac{\mu}{4}(t-t_{n_{0}})\right)\right)
\quad\quad
\forall t\in[t_{n_{0}},S).$$

Letting $t\to S^{-}$, we obtain that $\Ghat_{p}(S)<2c_{1}$, hence
$G_{p}(S)\leq 2\Ghat_{p}(S)<4c_{1}$ because of (\ref{th:Ghat-G}).
This contradicts the maximality of $S$ and completes the proof of
(\ref{hp:liminf-glob}).

\paragraph{\textmd{\textit{Faster decay of the range component}}}

Let us prove (\ref{th:range}).  To this end, we consider the
generalized Dirichlet quotients~(\ref{defn:Qd}) and (\ref{defn:Ghatd})
with $d:=2p$.  They are defined at least for every $t\geq t_{n_{0}}$
because of (\ref{th:u>0}).  Setting $d:=2p$ in (\ref{th:deriv-Ghat})
we obtain that
\begin{equation}
	\Ghat_{2p}'(t)\leq -\frac{\mu}{2}\Ghat_{2p}(t)+
	2(1+p)|u(t)|^{p}\cdot[G_{2p}(t)]^{1/2}\cdot\Ghat_{2p}(t)+
	\frac{2}{\delta}\frac{|g(t)|^{2}}{|u(t)|^{2+2p}}.	
	\label{est:Q2p'}
\end{equation}

From (\ref{hp:liminf-glob}) we deduce that
\begin{equation}
	|u(t)|^{p}\cdot[G_{2p}(t)]^{1/2}=
	|u(t)|^{p/2}\cdot[G_{p}(t)]^{1/2} \leq
	2\sqrt{c_{1}}|u(t)|^{p/2},
	\label{est:Q2p'-1}
\end{equation}
while from (\ref{th:Ghat-G}) and (\ref{hp:main-g-bis}) we deduce that
\begin{eqnarray}
	|g(t)|^{2} & \leq & 2K_{0}^{2}|u(t)|^{2+2p}\left(
	1+G_{2p}(t)\cdot|A^{1/2}u(t)|^{2q}\right)
	\nonumber  \\
	 & \leq & 2K_{0}^{2}|u(t)|^{2+2p}\left(
	1+2\Ghat_{2p}(t)\cdot|A^{1/2}u(t)|^{2q}\right).
	\label{est:Q2p'-2}
\end{eqnarray} 

Plugging (\ref{est:Q2p'-1}) and (\ref{est:Q2p'-2}) into 
(\ref{est:Q2p'}), we obtain that
$$\Ghat_{2p}'(t)\leq -\Ghat_{2p}(t)\cdot\left\{
\frac{\mu}{2}-c_{2}|u(t)|^{p/2}-
c_{3}|A^{1/2}u(t)|^{2q}
\right\}+c_{4}$$
for suitable constants $c_{2}$, $c_{3}$ and $c_{4}$.  Due to
assumption~(\ref{hp:u-limit}), there exists $T_{0}\geq t_{n_{0}}$ such
that
$$\frac{\mu}{2}-c_{2}|u(t)|^{p/2}-c_{3}|A^{1/2}u(t)|^{2q}
\geq\frac{\mu}{4} \quad\quad
\forall t\geq T_{0},$$
and hence
$$\Ghat_{2p}'(t)\leq -\frac{\mu}{4}\Ghat_{2p}(t)+c_{4}
\quad\quad
\forall t\geq T_{0}.$$

Integrating this differential inequality we conclude that
$\Ghat_{2p}(t)$ is uniformly bounded for every $t\geq T_{0}$.  Due to
(\ref{th:Ghat-G}), also $G_{2p}(t)$ is uniformly bounded for every
$t\geq T_{0}$, and this is enough to establish (\ref{th:range}).

\paragraph{\textmd{\textit{Slow decay of the solution}}}

Let us consider the function $y(t):=|u(t)|^{2}$. Since
$$|y'(t)|=2|\langle u'(t),u(t)\rangle|\leq
2|u'(t)|\cdot|u(t)|\leq
2G_{2p}(t)^{1/2}|u(t)|^{2+p},$$
from the uniform bound on $G_{2p}(t)$ we obtain that there exists a 
constant $c_{5}$ such that
$$|y'(t)|\leq c_{5}y(t)^{1+p/2}
\quad\quad
\forall t\geq T_{0},$$
and in particular
$$y'(t)\geq -c_{5}y(t)^{1+p/2}
\quad\quad
\forall t\geq T_{0}.$$

Since $y(T_{0})\neq 0$, integrating this differential inequality we
deduce (\ref{th:slow}).

\subsubsection*{Fast solutions}

In this last section of the proof it remains to consider the case
where (\ref{hp:u-neq-0}) holds true and $u(t)$ satisfies
(\ref{est:u<Au}) for some $T_{1}\geq 0$.  The constants $c_{6}$,
\ldots, $c_{16}$ which we introduce in the sequel are positive and
independent of time.

\paragraph{\textmd{\textit{Non-optimal exponential decay from below}}}

We prove that \begin{equation}
	|u'(t)|^{2}+|A^{1/2}u(t)|^{2}\geq c_{6}e^{-c_{7}t}
	\quad\quad
	\forall t\geq T_{1}.
	\label{th:u-exp-below}
\end{equation}

To this end, let us consider the usual hyperbolic energy $E(t)$
defined in (\ref{defn:E}).  Its time-derivative satisfies
\begin{equation}
	E'(t)=-4\delta|u'(t)|^{2}+2\langle g(t),u'(t)\rangle
	\geq -(4\delta+1)|u'(t)|^{2}-|g(t)|^{2}.
	\label{est:fast-F'}
\end{equation}

Let us estimate $|g(t)|$. Due to (\ref{est:u<Au}), inequality 
(\ref{hp:main-g-bis}) implies that
\begin{equation}
	|g(t)|^{2}\leq 
	2K_{0}^{2}\left(|u(t)|^{p}+|A^{1/2}u(t)|^{2q}\right)E(t)
	\quad\quad
	\forall t\geq T_{1}.
	\label{est:g-fast}
\end{equation}

Since $|u(t)|$ and $|A^{1/2}u(t)|$ are uniformly bounded because of
assumption (\ref{hp:u-limit}), it follows that $|g(t)|^{2}\leq
c_{8}E(t)$.  Plugging this estimate into (\ref{est:fast-F'}), we
deduce that
$$E'(t)\geq -c_{9}E(t)
\quad\quad
\forall t\geq T_{1}.$$

Integrating this differential inequality we obtain
(\ref{th:u-exp-below}).
 
\paragraph{\textmd{\textit{Non-optimal exponential decay from above}}}

We prove that there exists $T_{2}\geq T_{1}$ such that
\begin{equation}
	|u'(t)|^{2}+|A^{1/2}u(t)|^{2}\leq c_{10}e^{-c_{11}t}
	\quad\quad
	\forall t\geq T_{2}.
	\label{th:u-exp-above}
\end{equation}

To this end, let us consider the modified hyperbolic energy $\Ehat(t)$
defined in (\ref{defn:Fhat}). From (\ref{th:deriv-Fhat}), 
(\ref{est:g-fast}), and (\ref{th:Fhat-F})  it follows that
\begin{eqnarray*}
	\Ehat'(t) & \leq & -\frac{\mu}{2}\Ehat(t)+
	c_{12}\left(|u(t)|^{p}+|A^{1/2}u(t)|^{2q}\right)E(t) \\
	\noalign{\vspace{1ex}}
	 & \leq & -\left[\frac{\mu}{2}-2c_{12}
	\left(|u(t)|^{p}+|A^{1/2}u(t)|^{2q}\right)\right]\Ehat(t)
\end{eqnarray*}
for every $t\geq T_{1}$. Due to assumption (\ref{hp:u-limit}),
there exists $T_{2}\geq T_{1}$ such that
$$\frac{\mu}{2}-2c_{12}
\left(|u(t)|^{p}+|A^{1/2}u(t)|^{2q}\right)\geq
\frac{\mu}{4}
\quad\quad
\forall t\geq T_{2},$$
and hence
$$\Ehat'(t)\leq-\frac{\mu}{4}\Ehat(t)
\quad\quad
\forall t\geq T_{2}.$$

Integrating this differential inequality, and keeping 
(\ref{th:Fhat-F}) into account, we deduce that
$$E(t)\leq 2\Ehat(t)\leq 2\Ehat(T_{2})\exp
\left(-\frac{\mu}{4}(t-T_{2})\right)
\quad\quad
\forall t\geq T_{2},$$
which proves (\ref{th:u-exp-above}).

\paragraph{\textmd{\textit{Exact exponential decay rate}}}

To begin with, we observe that
\begin{equation}
	c_{13}e^{-c_{14}t}\leq
	|u'(t)|+|A^{1/2}u(t)|+|u(t)|\leq
	c_{15}e^{-c_{16}t}
	\quad\quad
	\forall t\geq T_{2}.
	\label{est:u-exp-bilateral}
\end{equation}

Indeed, the estimate from below is an immediate consequence of 
(\ref{th:u-exp-below}), while the estimate from above follows from 
(\ref{th:u-exp-above}) and (\ref{est:u<Au}). Now let us set
\begin{equation}
	r_{0}:=\sup\left\{r\geq 0:
	\lim_{t\to+\infty}
	\left(|u'(t)|+|u(t)|_{D(A^{1/2})}\right)
	e^{r t}=0\right\}.
	\label{defn:r0}
\end{equation}

From (\ref{est:u-exp-bilateral}) it follows that $r_{0}$ is finite and
strictly positive.  We claim that $r_{0}\in\mathcal{D}$, and that
there exists a nontrivial $r_{0}$-pure solution to the homogeneous
equation (\ref{eqn:lin-h}) for which (\ref{th:fast+}) holds true.

To begin with, we set
$$\beta_{1}:=\min\{(1+p)r_{0},(1+q)r_{0}\},$$ 
so that from assumption (\ref{hp:main-g}) we know now that
\begin{equation}
	\lim_{t\to +\infty}|g(t)|e^{\gamma t}=0
	\quad\quad
	\forall \gamma<\beta_{1}.
	\label{hp:s0}
\end{equation}

Let $\beta_{0}$ be defined by (\ref{defn:beta-0}), and let $\gamma$ be
any real number such that
\begin{equation}
	r_{0}<\gamma<\min\{\beta_{1},\beta_{0}\}.
	\label{hp:range-r}
\end{equation}

Since $\gamma<\beta_{1}$, from (\ref{hp:s0}) it follows that
$$|g(t)|\leq K_{g,\gamma}e^{-\gamma t}
\quad\quad
\forall t\geq 0$$
for a suitable constant $K_{g,\gamma}$.  Moreover, (\ref{hp:range-r})
implies that $\gamma\not\in\mathcal{D}$.  Therefore, we can apply
Proposition~\ref{prop:main-nh} with $\gamma_{0}:=\gamma$.  We deduce
that there exists a solution $w_{\gamma}(t)$ to the non-homogeneous
equation (\ref{eqn:lin-nh}) such that
\begin{equation}
	|w'_{\gamma}(t)|+|w_{\gamma}(t)|_{D(A^{1/2})}\leq 
	\Gamma_{0}K_{g,\gamma}e^{-\gamma t}
	\quad\quad
	\forall t\geq 0,
	\label{est:wr}
\end{equation}
and whose initial conditions satisfy suitable constraints.

We claim that $w_{\gamma}(t)$ does not depend on $\gamma$ as long as
(\ref{hp:range-r}) holds true.  This is almost trivial when
$r_{0}<\min\mathcal{D}$, because in this case
$\gamma<\beta_{0}=\min\mathcal{D}$ and hence $w_{\gamma}$ has initial
conditions $w_{\gamma}(0)=w_{\gamma}'(0)=0$.  If
$r_{0}\geq\min\mathcal{D}$, and $\alpha_{0}$ denotes the largest
element of $\mathcal{D}$ less than or equal to $r_{0}$, then
$w_{\gamma}(t)$ is uniquely characterized by the limit
(\ref{th:lin-nh-lim}) and by having $\alpha_{0}$-slow initial data,
and both conditions do not depend on $\gamma$ in the range
(\ref{hp:range-r}).

Therefore, in the sequel we denote $w_{\gamma}(t)$ just by $w(t)$ and, 
since (\ref{est:wr}) holds true for every $\gamma$ in the range 
(\ref{hp:range-r}), we deduce that
\begin{equation}
	\lim_{t\to +\infty}
	\left(|w'(t)|+|w(t)|_{D(A^{1/2})}\right)e^{\gamma t}=0
	\quad\quad
	\forall \gamma<\min\{\beta_{1},\beta_{0}\}.
	\label{est:lim-w}
\end{equation}

Now let us set $v(t):=u(t)-w(t)$.  Since $v(t)$ is a solution to the
homogeneous equation (\ref{eqn:lin-h}), it can be written as a finite
sum of $r$-pure solutions $v_{r}(t)$ to the same homogeneous equation,
with $r$ ranging over $\mathcal{D}$.  All terms $v_{r}(t)$ with
$r<r_{0}$ are necessarily equal to 0, because otherwise the supremum
in (\ref{defn:r0}) would be less than $r_{0}$.

We are now ready to prove our conclusions.  Let us assume by
contradiction that $r_{0}\not\in\mathcal{D}$.  In this case $v(t)$
is the sum of $r$-pure solutions $v_{r}(t)$ with $r>r_{0}$,
hence also $r\geq\beta_{0}$, and therefore
\begin{equation}
	\lim_{t\to +\infty}
	\left(|v'(t)|+|v(t)|_{D(A^{1/2})}\right)e^{\gamma t}=0
	\quad\quad
	\forall \gamma<\beta_{0}.
	\label{est:lim-v}
\end{equation}

But (\ref{est:lim-w}) and (\ref{est:lim-v}) imply that
$$\lim_{t\to +\infty}
\left(|u'(t)|+|u(t)|_{D(A^{1/2})}\right)e^{\gamma t}=0
\quad\quad
\forall \gamma<\min\{\beta_{1},\beta_{0}\},$$
so that the supremum in (\ref{defn:r0}) would be greater that 
$r_{0}$. This proves that $r_{0}\in\mathcal{D}$.  

Now $v(t)$ is the sum of an $r_{0}$-pure solution, which we denote by
$v_{0}(t)$, and possibly some other $r$-pure solutions $v_{r}(t)$ with
$r>r_{0}$, and hence $r\geq\beta_{0}$.  As a consequence,  it
turns out that
\begin{equation}
	\lim_{t\to +\infty}
	\left(|v'(t)-v_{0}'(t)|+|v(t)-v_{0}(t)|_{D(A^{1/2})}\right)
	e^{\gamma t}=0 \quad\quad
	\forall \gamma<\beta_{0}.
	\label{est:lim-v0}
\end{equation}

At this point (\ref{th:fast+}) follows from (\ref{est:lim-w}) and
(\ref{est:lim-v0}). Finally, $v_{0}(t)$ is not identically 0, because 
if not (\ref{est:lim-v}) would be true once again, and together with 
(\ref{est:lim-w}) this would contradict the maximality of $r_{0}$, 
exactly as before.\qed

\subsection{Proof of Theorem~\ref{thm:main-exponential}}

Let us first describe the plan of the proof, based on a fixed point
argument. Let us define $\beta_{0}$ as in (\ref{defn:beta-0}), and 
let us choose once for all a constant 
$s_{0}$ such that 
\begin{equation}
	s_{0}<r_{0}<(1+p)s_{0}<\min\{(1+p)r_{0},\beta_{0}\}.
	\label{defn:const}
\end{equation}

Let $v(t)$ be the $r_{0}$-pure solution to
the homogeneous equation (\ref{eqn:lin-h}) with initial data
$(v_{0},v_{1})$.  Let $z(t)$ be the $r_{0}$-fast solution to
the homogeneous equation (\ref{eqn:lin-h}) with initial data
$(z_{0},z_{1})$. Let
$K_{1}$ and $K_{2}$ be two constants such that
\begin{equation}
	|v(t)|_{D(A^{1/2})}\leq 
	K_{1}\left(|v_{1}|+|v_{0}|_{D(A^{1/2})}\right)e^{-s_{0}t}
	\quad\quad
	\forall t\geq 0,
	\label{est:fast-v}
\end{equation}
\begin{equation}
	|z(t)|_{D(A^{1/2})}\leq 
	K_{2}\left(|z_{1}|+|z_{0}|_{D(A^{1/2})}\right)e^{-s_{0}t}
	\quad\quad
	\forall t\geq 0.
	\label{est:fast-z}
\end{equation}

The constants $K_{1}$ and $K_{2}$ exist because the left-hand side of
(\ref{est:fast-v}) decays at least as $(1+t)e^{-r_{0}t}$, and the
left-hand side of (\ref{est:fast-z}) decays at least as
$(1+t)e^{-\beta_{0}t}$ (or it is identically 0 if $\beta_{0}=+\infty$).  

Let $\Gamma_{0}$ be the constant which appears in
Proposition~\ref{prop:main-nh} when $\gamma_{0}:=(1+p)s_{0}$, and let
us assume that $\ep_{0}>0$ is small enough so that
\begin{equation}
	(K_{1}+K_{2}+1)\ep_{0}<R_{0},
	\label{hp:small-1}
\end{equation}
\begin{equation}
	2L\Gamma_{0}(K_{1}+K_{2}+1)^{1+p}\ep_{0}^{p}<1.
	\label{hp:small-2}
\end{equation}

Let us consider the space 
$$\mathbb{X}:=\left\{\psi\in C^{0}\left([0,+\infty);D(A^{1/2})\right):
|\psi(t)|_{D(A^{1/2})}\leq\ep_{0}\quad\forall t\geq 0\right\}.$$

It is well-known that $\mathbb{X}$ is a complete metric space with respect 
to the distance
$$\mbox{dist}(\psi_{1},\psi_{2}):=\sup
\left\{|\psi_{1}(t)-\psi_{2}(t)|_{D(A^{1/2})}:t\geq 0\right\}.$$

For every $\psi\in \mathbb{X}$ we set
$$g_{\psi}(t):=f\left(v(t)+z(t)+\psi(t)e^{-\gamma_{0}t}\right)
\quad\quad
\forall t\geq 0,$$
and we consider the non-homogeneous linear equation
\begin{equation}
	w''(t)+2\delta w'(t)+Aw(t)=g_{\psi}(t).
	\label{eqn:wg}
\end{equation}

We claim that this equation admits a unique solution $w_{\psi}(t)$ such
that
\begin{equation}
	\lim_{t\to +\infty}
	\left(|w_{\psi}'(t)|+|w_{\psi}(t)|_{D(A^{1/2})}\right)
	e^{r_{0}t}=0,
	\label{th:wg-limit}
\end{equation}
and whose initial data are $r_{0}$-slow. Finally, we set
\begin{equation}
	\overline{\psi}(t):=w_{\psi}(t)e^{\gamma_{0} t}
	\quad\quad
	\forall t\geq 0,
	\label{defn:psi-bar}
\end{equation}
and we claim that the following three statements hold true, provided
that the smallness assumptions (\ref{hp:small-1}) and
(\ref{hp:small-2}) are satisfied.

\begin{itemize}
	\item \emph{Well-posedness of the construction}.  The functions
	$g_{\psi}$ and $w_{\psi}$ are well-defined for every $\psi\in
	\mathbb{X}$.  
	
	\item \emph{Closedness}.  It turns out that $\overline{\psi}\in \mathbb{X}$
	for every $\psi\in \mathbb{X}$.

	\item \emph{Contractivity}.  The map $\mathcal{F}:\mathbb{X}\to
	\mathbb{X}$ defined by $\mathcal{F}(\psi)=\overline{\psi}$ is a
	contraction.
\end{itemize}

If we prove the three claims, then the conclusion follows.
Indeed the contractivity implies that $\mathcal{F}$ has a fixed point,
namely there exists $\psi\in\mathbb{X}$ such that $\overline{\psi}=\psi$. If 
$\psi$ is the fixed point, then the function defined by
$$u_{\psi}(t):=v(t)+z(t)+w_{\psi}(t)$$
is the solution we were looking for.

Indeed, (\ref{th:limit}) holds true because both $z(t)$ and 
$w_{\psi}(t)$ decay faster than $e^{-r_{0}t}$. The initial conditions 
of $u_{\psi}(t)$ are of the form (\ref{th:u-data}), where 
$(w_{0},w_{1})$ are the initial data of $w_{\psi}(t)$, which are 
$r_{0}$-slow. Finally, since $\overline{\psi}(t)=\psi(t)$, from 
(\ref{defn:psi-bar}) it follows that
$$v(t)+z(t)+\psi(t)e^{-\gamma_{0}t}=
v(t)+z(t)+\overline{\psi}(t)e^{-\gamma_{0}t}=
v(t)+z(t)+w_{\psi}(t)=u_{\psi}(t),$$
and hence 
$$g_{\psi}(t)=f\left(v(t)+z(t)+\psi(t)e^{-\gamma_{0}t}\right)=
f(u_{\psi}(t)).$$

Since $v(t)$ and $z(t)$ are solutions to the homogeneous equation, we 
conclude that
$$u_{\psi}''(t)+2\delta u_{\psi}'(t)+Au_{\psi}(t)=
w_{\psi}''(t)+2\delta w_{\psi}'(t)+Aw_{\psi}(t)=
g_{\psi}(t)=f(u_{\psi}(t)),$$
which proves that $u_{\psi}(t)$ is a solution to (\ref{eqn:semilin}).

\subparagraph{\textmd{\textit{Well-posedness of the construction}}}

From (\ref{est:fast-v}), (\ref{est:fast-z}), and our definition of
$\mathbb{X}$, it turns out that 
\begin{eqnarray*}
	\left|v(t)+z(t)+\psi(t)e^{-\gamma_{0}t}\right|_{D(A^{1/2})} & \leq &
	K_{1}\left(|v_{1}|+|v_{0}|_{D(A^{1/2})}\right)e^{-s_{0}t} \\
	 & & \mbox{}+
	 K_{2}\left(|z_{1}|+|z_{0}|_{D(A^{1/2})}\right)e^{-s_{0}t}+
	 \ep_{0}e^{-\gamma_{0}t}.
\end{eqnarray*}

Since $s_{0}<\gamma_{0}$, from (\ref{hp:fast}) we obtain that
\begin{equation}
	\left|v(t)+z(t)+\psi(t)e^{-\gamma_{0}t}\right|_{D(A^{1/2})}\leq
	(K_{1}+K_{2}+1)\ep_{0}e^{-s_{0}t}
	\quad\quad
	\forall t\geq 0.
	\label{est:arg-phig}
\end{equation}

Due to the smallness assumption (\ref{hp:small-1}), this means in 
particular that
$$\left|v(t)+z(t)+\psi(t)e^{-\gamma_{0}t}\right|_{D(A^{1/2})}\leq
(K_{1}+K_{2}+1)\ep_{0}<R_{0},$$
which proves that $g_{\psi}(t)$ is well-defined.

Setting $v=0$ into~(\ref{hp:f-plip}), from (\ref{hp:f0}) we obtain
that
$$|f(u)|\leq L|u|^{1+p}_{D(A^{1/2})}
\quad\quad
\forall u\in B_{R_{0}}.$$

Therefore, from (\ref{est:arg-phig}) we deduce that
\begin{equation}
	|g_{\psi}(t)|\leq L(K_{1}+K_{2}+1)^{1+p}
	\ep_{0}^{1+p}e^{-(1+p)s_{0}t}.
	\label{est:phig}
\end{equation}

Now we apply Proposition~\ref{prop:main-nh} with
$\gamma_{0}:=(1+p)s_{0}$.  Due to (\ref{defn:const}), we are in the
case where $\gamma_{0}>\min\mathcal{D}$, and $r_{0}$ is the largest
element of $\mathcal{D}$ smaller than $\gamma_{0}$ (namely what was
called $\alpha_{0}$ in Proposition~\ref{prop:main-nh}).  We obtain
that (\ref{eqn:wg}) has a unique solution $w_{\psi}(t)$ with
$r_{0}$-slow initial data and such that (\ref{th:wg-limit}) holds
true.  This proves that $w_{\psi}(t)$ is well-defined.

\subparagraph{\textmd{\textit{Closedness}}}

To begin with, we observe that $\overline{\psi}:[0,+\infty)\to
D(A^{1/2})$ is a continuous map because of the regularity of
$w_{\psi}$ and $f$.  Due to estimate (\ref{est:phig}) and the
smallness assumption (\ref{hp:small-2}), from
Proposition~\ref{prop:main-nh} we obtain that
$$|w_{\psi}(t)|_{D(A^{1/2})}\leq \Gamma_{0}L(K_{1}+K_{2}+1)^{1+p}
\ep_{0}^{1+p}e^{-(1+p)s_{0}t}\leq\ep_{0}e^{-\gamma_{0}t},$$
from which we conclude that 
$|\overline{\psi}(t)|_{D(A^{1/2})}\leq\ep_{0}$ for every $t\geq 0$. This 
proves that $\overline{\psi}\in\mathbb{X}$.

\subparagraph{\textmd{\textit{Contractivity}}}

Let $\psi_{1}$ and $\psi_{2}$ be two elements of $\mathbb{X}$.
Estimate (\ref{est:arg-phig}) holds true also with $\psi_{1}$ and
$\psi_{2}$ instead of $\psi$.  Therefore, from (\ref{hp:f-plip}) we
deduce that
\begin{eqnarray*}
	\left|g_{\psi_{1}}(t)-g_{\psi_{2}}(t)\right|_{D(A^{1/2})} & \leq
	& 2L(K_{1}+K_{2}+1)^{p}\ep_{0}^{p}e^{-ps_{0}t}
	\left|\psi_{1}(t)-\psi_{2}(t)\right|_{D(A^{1/2})}e^{-\gamma_{0}t}
	\\
	 & \leq & 2L(K_{1}+K_{2}+1)^{p}\ep_{0}^{p}e^{-\gamma_{0}t}
	 \cdot\mbox{dist}(\psi_{1},\psi_{2}).
\end{eqnarray*}

Let $w_{\psi_{1}}(t)$ and $w_{\psi_{2}}(t)$ denote the corresponding
solutions to (\ref{eqn:wg}) in the sense of
Proposition~\ref{prop:main-nh}.  Since
$w_{\psi_{1}}(t)-w_{\psi_{2}}(t)$ solves the same equation with
forcing term $g_{\psi_{1}}(t)-g_{\psi_{2}}(t)$, now (\ref{th:lin-nh})
reads as
$$\left|w_{\psi_{1}}(t)-w_{\psi_{2}}(t)\right|_{D(A^{1/2})}\leq
\Gamma_{0}\cdot 2L(K_{1}+K_{2}+1)^{p} \ep_{0}^{p}\cdot
e^{-\gamma_{0}t}\cdot \mbox{dist}(\psi_{1},\psi_{2}).$$
	
Multiplying by $e^{\gamma_{0}t}$, and taking the supremum for $t\geq 
0$, we finally obtain that
$$\mbox{dist}(\overline{\psi}_{1},\overline{\psi}_{2})\leq
2\Gamma_{0}L(K_{1}+K_{2}+1)^{p}\ep_{0}^{p}
\cdot\mbox{dist}(\psi_{1},\psi_{2}),$$
so that the smallness assumption (\ref{hp:small-2}) implies that the
map $\mathcal{F}:\mathbb{X}\to \mathbb{X}$ is a contraction.\qed

\setcounter{equation}{0}
\section{Applications}\label{sec:applications}

\subsection{Semilinear dissipative hyperbolic equations}

Let $\Omega\subseteq\re^{n}$ be a bounded connected open set with
Lipschitz boundary (or any other condition which guarantees Sobolev
embeddings).  As a model case, we consider dissipative hyperbolic
equations of the form
\begin{equation}
	u_{tt}+2\delta u_{t}-\Delta u\pm |u|^{p}u=0
	\quad\quad
	\mbox{in }\Omega\times [0,+\infty),
	\label{pbm:n-psi}
\end{equation}
with homogeneous Neumann boundary conditions, or of the form
\begin{equation}
	u_{tt}+2\delta u_{t}-\Delta u-\lambda u\pm |u|^{p}u=0
	\quad\quad
	\mbox{in }\Omega\times [0,+\infty),
	\label{pbm:d-psi}
\end{equation}
with homogeneous Dirichlet boundary conditions.  In both cases, $|u|$
denotes the absolute value of $u$.  In the case of equation
(\ref{pbm:d-psi}), we assume that $\lambda\leq\lambda_{1}(\Omega)$,
where $\lambda_{1}(\Omega)$ denotes the first eigenvalue of $-\Delta$
with homogeneous Dirichlet boundary conditions in $\Omega$.

The functional setting is the classical one, namely $H:=L^{2}(\Omega)$
and $Au:=-\Delta u$ with a suitable domain depending on boundary
conditions.  We refer to~\cite{ggh:sol-lentes} or~\cite{ggh:casc-par}
for further details.  We just point out that $A$ is a coercive
operator in the subcritical Dirichlet case where
$\lambda<\lambda_{1}(\Omega)$, but it is just a nonnegative operator
both in the Neumann case (where the kernel of $A$ is the space of
constant functions), and in the critical Dirichlet case where
$\lambda=\lambda_{1}(\Omega)$ (where the kernel of $A$ is the
eigenspace of $-\Delta$ relative to the eigenvalue
$\lambda_{1}(\Omega)$).

As for the nonlinear term, we set
$$[f(u)](x):=\mp|u(x)|^{p}u(x)
\quad\quad
\forall x\in\Omega.$$

The function $f$ satisfies the assumptions of our abstract results
provided that the Sobolev embedding $H^{1}(\Omega)\subseteq
L^{2+2p}(\Omega)$ holds true.  In turn, this condition is satisfied
for every $p>0$ if $n\in\{1,2\}$, and when $0<p\leq 2/(n-2)$ if $n\geq
3$.  We refer to Section~4.1 of~\cite{ggh:sol-lentes} for a proof of
these basic facts, which are independent of the sign in
(\ref{pbm:n-psi}) or (\ref{pbm:d-psi}).

The sign becomes relevant when looking for global solutions for all 
initial data in the energy space. With the ``right sign'', from the 
results of~\cite{ggh:sol-lentes} and of the present paper, we can 
prove the following statement.

\begin{thm}[Right sign]\label{thm:right-sign}
	Let $\Omega\subseteq\re^{n}$ be a bounded open set with Lipschitz
	boundary, and let $p$ be a positive exponent, with no further
	restriction if $n\in\{1,2\}$, and $p\leq 2/(n-2)$ if $n\geq 3$.
	
	Then the following five statements apply to both the Neumann
	problem for equation~(\ref{pbm:n-psi}) with the plus sign, and to
	the Dirichlet problem for equation~(\ref{pbm:d-psi}) with
	the plus sign and $\lambda=\lambda_{1}(\Omega)$.  It is intended
	that $D(A^{1/2})=H^{1}(\Omega)$ in the Neumann case, and
	$D(A^{1/2})=H^{1}_{0}(\Omega)$ in the Dirichlet case.
	\begin{enumerate}
		\renewcommand{\labelenumi}{(\arabic{enumi})}
		
		\item \emph{(Global existence and uniqueness)} For every
		$(u_{0},u_{1})\in D(A^{1/2})\times H$ there exists a unique
		global solution with the regularity (\ref{th:u-reg}).
		
		\item \emph{(Decay estimate from above)} All solutions satisfy
		\begin{equation}
			\|u(t)\|_{L^{2}(\Omega)}\leq \frac{C}{(1+t)^{1/p}}
			\quad\quad
			\forall t\geq 0			
			\label{th:rs-u}
		\end{equation}
		and
		\begin{equation}
			\|u'(t)\|_{L^{2}(\Omega)}+
			\|\nabla u(t)\|_{L^{2}(\Omega)}\leq 
			\frac{C}{(1+t)^{1+1/p}}
			\quad\quad
			\forall t\geq 0			
			\label{th:rs-u'}
		\end{equation}
		for a suitable constant $C$ (depending on the solution).
		
		\item \emph{(Classification of decay rates)} 
		All non-zero solutions are either slow or fast in the sense 
		of Theorem~\ref{thm:main-alternative}.
		
		\item \emph{(Existence of slow solutions)} There exists a
		nonempty open set $\mathcal{S}\subseteq D(A^{1/2})\times H$
		such that all solutions with initial data in $S$ are slow.
		
		\item \emph{(Existence of fast solutions)} There exists
		families of fast solutions parametrized in the sense of
		Theorem~\ref{thm:main-exponential}.
		
	\end{enumerate}
		
\end{thm}

Of course the theory applies also to the Dirichlet problem with 
$\lambda<\lambda_{1}(\Omega)$, but in that case the operator is 
coercive and we have only fast solutions.

Let us spend a few words on the proof of Theorem~\ref{thm:right-sign}.
Statement~(1) is a well-known result.  Estimate (\ref{th:rs-u}) was
proved in Theorem~2.2 of~\cite{ggh:sol-lentes} together with a weaker
version of (\ref{th:rs-u'}), in which the exponent $(1+1/p)$ is
replaced by $(1/2+1/p)$.  This weaker estimate is enough to conclude
that $u(t)$ decays to 0 in $D(A^{1/2})$, which allows to apply
Theorem~\ref{thm:main-alternative} of the present paper to all
solutions.  Thus we obtain statement~(3), and also estimate
(\ref{th:rs-u'}) with the correct exponent, which follows from
(\ref{th:range}) in the case of slow solutions, and is trivially true
both for the null solution and for fast solutions, which decay
exponentially.  Statement~(4) is a consequence of Theorem~2.3 
of~\cite{ggh:sol-lentes}. Statement~(5) follows from 
Theorem~\ref{thm:main-exponential} of the present paper.

When the nonlinear term has the wrong sign, global existence is known
only in special cases, for example when the origin falls in the
so-called potential well.  This technique requires that the operator
is coercive and controls the nonlinear term (which means Sobolev
embeddings).  Since the coerciveness of the operator is essential,
this theory applies neither to the Neumann case, nor to the critical
Dirichlet case.  In other words, the potential well applies only to
the subcritical Dirichlet case, in which case we obtain the following
result.

\begin{thm}[Wrong sign, with potential well]
	Let $\Omega\subseteq\re^{n}$ be a bounded open set with Lipschitz
	boundary.  Let $p$ be a positive exponent, with no further
	restriction if $n\in\{1,2\}$, and $p\leq 2/(n-2)$ if $n\geq 3$.
	Let us consider the Dirichlet problem for equation
	(\ref{pbm:d-psi}) with the minus sign and
	$\lambda<\lambda_{1}(\Omega)$.
	
	Then there exists $R_{0}>0$ such that, for every $u_{0}\in
	B_{R_{0}}$ (defined as in (\ref{defn:ball})), the problem has a
	unique global solution with the regularity (\ref{th:u-reg}).
	
	Moreover, every non-zero solution in $B_{R_{0}}$ is fast in the
	sense of Theorem~\ref{thm:main-alternative}, and there exist
	families of fast solutions parametrized in the sense of
	Theorem~\ref{thm:main-exponential}.
\end{thm}

When there is no potential well, Theorem~\ref{thm:main-alternative}
keeps on classifying all possible decay rates of those solutions which
exist globally and decay.  On the other hand, nothing in this case
guarantees decay, or even global existence, of solutions.

Nevertheless, there is one notable exception.
Theorem~\ref{thm:main-exponential} provides families of global
solutions with exponential decay without assuming neither the
coercivity of the operator, nor sign conditions on the nonlinear term.
Thus we obtain the following existence result.  

\begin{thm}[Wrong sign, without potential well]\label{thm:wrong-sign}
	Let $\Omega\subseteq\re^{n}$ be a bounded open set with Lipschitz
	boundary.  Let $p$ be a positive exponent, with no further
	restriction if $n\in\{1,2\}$, and $p\leq 2/(n-2)$ if $n\geq 3$.
	
	Let us consider the Neumann problem for equation~(\ref{pbm:n-psi})
	or the Dirichlet problem for equation~(\ref{pbm:d-psi}) with
	$\lambda=\lambda_{1}(\Omega)$.
	
	Then there exist families of fast solutions parametrized in the
	sense of Theorem~\ref{thm:main-exponential}, independently of the
	sign in the nonlinear term.
		
\end{thm}

We conclude by pointing out that our abstract results apply also to 
equations with second order operators with non-constant coefficients, 
or with higher order operators such as $\Delta^{2}$. We also allow 
more general nonlinear terms depending on $x$ and $t$.

\subsection{Degenerate Kirchhoff equations in finite dimension}

In this final section we present a different application of our 
theory. We consider a degenerate Kirchhoff equation
\begin{equation}
	u''(t)+u'(t)+\left|B^{1/2}u(t)\right|^{2\alpha}Bu(t)=0,
	\label{eqn:K}
\end{equation}
where $\alpha$ is a positive real number and $B$ is a self-adjoint
operator on a Hilbert space $H$.  Equations of this type have long
been considered in the literature, but only partial results are known.
As for global existence, the main result is that a global solution
exists provided that initial data $(u_{0},u_{1})\in D(B)\times
D(B^{1/2})$ satisfy the nondegeneracy condition $B^{1/2}u_{0}\neq 0$,
and a suitable smallness assumption.  This was proved in~\cite{ny} in
the case $\alpha\geq 1$, and then in~\cite{ghisi:JDE2006} in the case
$0<\alpha<1$.

As for decay estimates, let us assume that the operator $B$ is
coercive, because if not solutions do not necessarily decay to 0 (just
think to the limit case where $B$ is the null operator).  Under this
coerciveness assumption, it is know that solutions provided in
literature satisfy
\begin{equation}
	\frac{C_{1}}{(1+t)^{1/(2\alpha)}}\leq|u(t)|\leq
	\frac{C_{2}}{(1+t)^{1/(2\alpha)}}
	\label{est:k-slow}
\end{equation}
for suitable positive constants $C_{1}$ and $C_{2}$, which means that
these solutions are slow.  Analogous estimates hold true for
$|B^{1/2}u(t)|$ and $|Bu(t)|$.  This was proved
in~\cite{ny,gg:JDE2008,ghisi:SIAM}.

On the other hand, there exist solutions to (\ref{eqn:K}) which are
not slow.  For example, if we limit ourselves to simple modes, namely
solutions of the form $u(t):=u_{k}(t)e_{k}$, where $e_{k}$ is an
eigenvector of $B$ corresponding to a positive eigenvalue
$\lambda_{k}$, then (\ref{eqn:K}) reduces to the ordinary differential
equation $$u_{k}''(t)+u_{k}'(t)+
\lambda_{k}^{\alpha+1}|u_{k}(t)|^{2\alpha}u_{k}(t)=0,$$
and it is well-known after~\cite{h:ode} that this equation admits both
slow solutions decaying as $t^{-1/(2\alpha)}$ and fast solutions
decaying as $e^{-t}$.

Now we can say that this alternative holds true more generally for 
solutions with a finite number of modes, or more generally for 
solutions living in a subspace of $H$ where (the restriction of) $B$ 
is a bounded operator. 

\begin{thm}
	Let $H$ be a separable Hilbert space, and let $B$ be a linear
	operator on $H$.  Let us assume that $B$ is bounded, symmetric and
	coercive.
	
	Then for every $\alpha> 0$ the following 
	conclusions hold true.
	\begin{enumerate}
		\renewcommand{\labelenumi}{(\arabic{enumi})}
		
		\item \emph{(Global existence and uniqueness)} For every
		$(u_{0},u_{1})\in H\times H$, problem
		(\ref{eqn:K})--(\ref{eqn:u-data}) admits a unique global
		solution $u\in C^{2}([0,+\infty),H)$.
		
		\item \emph{(Classification of decay rates)} Every non-zero
		solution is either a slow solution satisfying
		(\ref{est:k-slow}) for suitable positive constants $C_{1}$ 
		and $C_{2}$, or a fast solution for which there exists 
		$v_{0}\in H$, with $v_{0}\neq 0$, such that
		\begin{equation}
			\lim_{t\to +\infty}
			\left(|u'(t)+v_{0}e^{-t}|+|u(t)-v_{0}e^{-t}|\right)
			e^{\gamma t}=0
			\quad\quad
			\forall\gamma<1+2\alpha.
			\label{th:k-fast}
		\end{equation}
		
		\item \emph{(Existence of slow solutions)} There exists a 
		nonempty open set $\mathcal{S}\subseteq H\times H$ of initial 
		data originating slow solutions.
		
		\item \emph{(Existence of fast solutions)} For every 
		$v_{0}\in H$, small enough but different from 0, there exists 
		at least one solution satisfying (\ref{th:k-fast}).
		
	\end{enumerate}
\end{thm}

Let us sketch the proof, which is just an application of our theory.
Statement~(1) follows from the boundedness of $B$ and the fact that 
the usual Hamiltonian
$$|u'(t)|^{2}+\frac{1}{\alpha+1}|B^{1/2}u(t)|^{2(\alpha+1)}$$
is constant along trajectories. 

Now let us rewrite (\ref{eqn:K}) in the form
$$u''(t)+u'(t)=-\left|B^{1/2}u(t)\right|^{2\alpha}Bu(t)
=:f(u(t)).$$

It can be seen that $f(u)=-\nabla F(u)$ with
$F(u):=(\alpha+1)^{-1}|B^{1/2}u|^{2(\alpha+1)}$, and that $|f(u)|\leq
K_{0}|u|^{2\alpha+1}$ for a suitable constant $K_{0}$ because the
norms in $H$, $D(B^{1/2})$ or $D(B)$ are equivalent due to the
boundedness and coerciveness of $B$.  Therefore, equation
(\ref{eqn:K}) fits in the abstract framework of \cite{ggh:sol-lentes}
and of the present paper with $A$ equal to the null operator and
$p:=2\alpha$.  At this point, from Theorem~2.2
of~\cite{ggh:sol-lentes} it follows that all solutions satisfy
$$|u(t)|\leq\frac{C_{3}}{(1+t)^{1/(2\alpha)}} 
\quad\quad
\forall t\geq 0$$
for a suitable constant $C_{3}$.  In particular, all solutions decay
to 0, and hence we can apply Theorem~\ref{thm:main-alternative} of the
present paper, which gives the slow-fast alternative.  As for the
asymptotic profile of fast solutions, it is enough to remark that now
the associated homogeneous equation is $u''+u'=0$, so that the only
positive element of $\mathcal{D}$ is $r_{0}=1$, and $r_{0}$-pure
solutions are of the form $v_{0}(t)=v_{0}e^{-t}$ for some $v_{0}\in
H$.  This proves statement~(2).

Statement~(3) follows from Theorem~2.3 of~\cite{ggh:sol-lentes}.

Statement~(4) follows from Theorem~\ref{thm:main-exponential} of the 
present paper applied with $r_{0}=1$, after observing the structure 
of $r_{0}$-pure solutions and the fact that the only $r_{0}$-fast 
initial datum is $(0,0)$. 

\subsubsection*{\centering Acknowledgments}

This work has been done while the first two authors were visiting the
Laboratoire Jacques Louis Lions of the UPMC (Paris~VI).  The stay was
partially supported by the FSMP (Fondation Sciences Math\'{e}matiques
de Paris).  The first two authors are members of the Gruppo Nazionale
per l'Analisi Matematica, la Probabilit\`{a} e le loro Applicazioni
(GNAMPA) of the Istituto Nazionale di Alta Matematica (INdAM).

\label{NumeroPagine}


{\small \begin{thebibliography}{99}

	\bibitem{ghisi:JDE2006}{\sc M.\ Ghisi}; Global solutions for dissipative
	Kirchhoff strings with non-Lipschitz nonlinear term.  \emph{J.\
	Differential Equations} \textbf{230} (2006), no.~1, 128--139.

	\bibitem{ghisi:SIAM}{\sc M.\ Ghisi}; Asymptotic limits for mildly
	degenerate Kirchhoff equations.  \emph{SIAM J.\ Math.\
	Anal.}\ \textbf{45} (2013), no.~3, 1886--1906.

	\bibitem{gg:JDE2008}{\sc M.\ Ghisi, M.\ Gobbino}; Hyperbolic-parabolic
	singular perturbation for mildly degenerate Kirchhoff equations:
	time-decay estimates.  \emph{J.\ Differential Equations} \textbf{245}
	(2008), no.~10, 2979--3007.

	\bibitem{ggh:sol-lentes}{\sc M.\ Ghisi, M.\ Gobbino, A.~Haraux};
	Optimal decay estimates for the general solution to a class of
	semi-linear dissipative hyperbolic equations.  To appear on
	\emph{J.\ Eur.\ Math.\ Soc.\ (JEMS)}.  Preprint
	\texttt{arXiv:1306.3644}.

	\bibitem{ggh:casc-par}{\sc M.\ Ghisi, M.\ Gobbino, A.~Haraux};
	A description of all possible decay rates for solutions of 
	some semilinear parabolic equations.  To appear on
	\emph{J.\ Math.\ Pures Appl.}\ Preprint
	\texttt{arXiv:1402.5355}.
	
	\bibitem{h:hyp}{\sc A.\ Haraux, M.\ A.\ Jendoubi}; Decay estimates
	to equilibrium for some evolution equations with an analytic
	nonlinearity, \emph{ Asymptot.\ Anal.}\ \textbf{26} (2001), no.~1,
	21--36.

	\bibitem{h:ode}{\sc A.\ Haraux}; Slow and fast decay of
	solutions to some second order evolution equations.  \emph{J.\
	Anal.\ Math.}\ \textbf{95} (2005), 297--321.
	
	\bibitem{h:nodea} \textsc{A.~Haraux}; Decay rate of the range
	component of solutions to some semilinear evolution equations.
	\emph{NoDEA Nonlinear Differential Equations Appl.}\ \textbf{13}
	(2006), no.~4, 435--445.
	
	\bibitem{ny}{\sc K.\ Nishihara, Y.\ Yamada}; On global solutions
	of some degenerate quasilinear hyperbolic equations with
	dissipative terms.  \emph{Funkcial.\ Ekvac.}\ \textbf{33} (1990),
	no.~1, 151--159.

	\bibitem{hartmann-1} \textsc{P.~Hartman}; On local homeomorphisms
	of Euclidean spaces.  \emph{Bol.\ Soc.\ Mat.\ Mexicana (2)}
	\textbf{5} (1960), 220--241.

	\bibitem{hartmann-2} \textsc{P.~Hartman}; A lemma in the theory of
	structural stability of differential equations.  \emph{Proc.\
	Amer.\ Math.\ Soc.}\ \textbf{11} (1960), 610--620.

	\bibitem{hartmann-3} \textsc{H.~M.~Rodrigues, J.~Sol\`{a}-Morales};
	Linearization of class $C^{1}$ for contractions on Banach spaces.
	\emph{J.\ Differential Equations} \textbf{201} (2004), no.~2,
	351--382.

\end{thebibliography}
}
\end{document}